\documentclass[11pt]{amsart}
\usepackage{latexsym}
\usepackage{amsmath,amsthm,amssymb}
\usepackage[dvips]{graphicx}
\usepackage{tabularx}
\newcommand{\pf}{\textit{Proof.} \ }
\newcommand{\remark}{\textit{Remark.} \ }
\newcommand{\example}{\textit{Example.} \ }
\theoremstyle{plain}
\newtheorem{df}{Definition}[section]
\newtheorem{thm}[df]{Theorem}

\newtheorem{prop}[df]{Proposition}
\newtheorem{lem}[df]{Lemma}
\newtheorem{cor}[df]{Corollary}

\def\dfrac#1#2{{\displaystyle\frac{#1}{#2}}}

\newcommand{\R}{\boldsymbol{R}}

\newcommand{\C}{\boldsymbol{C}}

\newcommand{\Adj}{\mathrm{Ad}}
\def\l{{\mathfrak{l}}}
\newcommand{\h}{\mathfrak{h}}

\newcommand{\g}{\mathfrak{g}}

\newcommand{\da}{\mathfrak{a}}
\newcommand{\db}{\mathfrak{b}}
\newcommand{\sll}{\mathfrak{sl}}
\newcommand{\gl}{\mathfrak{gl}}

\usepackage{enumerate}
\usepackage{multirow,bigdelim}
\usepackage{delarray}

\date{}
\begin{document}
\title{Left invariant flat projective structures on Lie groups and prehomogeneous vector spaces}
\author{ Hironao Kato}
\address{Department of Mathematics, Graduate School of Science, Hiroshima University, Higashi-Hiroshima 739-8526, Japan}
\email{katoh-in-math@hiroshima-u.ac.jp}
\begin{abstract}
We show the correspondence between left invariant flat projective structures on Lie groups and certain prehomogeneous vector spaces. 
Moreover by using the classification theory of prehomogeneous vector spaces, we classify complex Lie groups admitting irreducible left 
invariant flat complex projective structures. As a result, direct sums of special 
linear Lie algebras $\sll(2) \oplus \sll(m_1) \oplus \cdots \oplus \sll(m_k)$ admit left invariant flat complex projective structures if 
the equality $4 + m_1^2 + \cdots + m_k^2 -k - 4 m_1 m_2 \cdots m_k = 0$ holds.  
These contain $\sll(2)$, $\sll(2) \oplus \sll(3)$, $\sll(2) \oplus \sll(3) \oplus \sll(11)$ for example. 
\end{abstract}
\maketitle
\section{Introduction}
A flat projective structure on a manifold is a maximal atlas whose charts take values in the projective space and 
coordinate changes are projective transformations (cf. \cite{goldman}). Our definition exactly agrees with the familiar 
one using projective equivalence classes of connections, which is explained in \S 2 of this paper. 
A flat projective structure on a Lie group is said to be left invariant if coordinate 
expressions of left translations are projective transformations (see $\S 2$). We abbreviate a left invariant flat projective structure to 
IFPS. Likewise on complex Lie groups we can consider left invariant flat complex projective structures (abbrev. complex IFPS) by taking the 
complex projective space as a model space. 
Then there arises a natural problem; on a given (complex) Lie group, is there a (complex) IFPS or not? 
Concerning this problem, Agaoka \cite{agaoka}, Urakawa \cite{urakawa}, Elduque \cite{elduque} proved that a real simple Lie group admits an 
IFPS if and only if its Lie algebra is $\sll(n+1, \R)$ or $\mathfrak{su}^*(2n)$ ($n \geq 1$).   
However concerning real and complex semisimple Lie groups, the classification problem is open.  

In this paper by using the theory of prehomogeneous vector spaces, 
we give an answer to this classification problem under one restriction on geometric structures called irreducibility 
(cf. Definition \ref{irreducibility}). 
Note that IFPSs are divided into two groups, that is, irreducible IFPSs and reducible ones.   
Our main theorem is stated in the following form: 
\begin{thm}\label{main thm}
A complex Lie group admits an irreducible complex IFPS if and only if 
its Lie algebra is of the form 
\[\sll(a) \oplus \sll(m_1) \oplus \cdots \oplus \sll(m_k)\]  
where $a = 2$, $3$, or $5$ \ $(k \geq 1, \ m_i \geq 1)$ and satisfies the equality 
\[(\ast \ast) \quad a^2 + m_1^2 + \cdots + m_k^2 -k - 2 a m_1 m_2 \cdots m_k = 0.\]  
\end{thm}
These Lie algebras include an infinite number of Lie algebras: $\sll(2)$, $\sll(2) \oplus \sll(3)$, 
$\sll(2) \oplus \sll(3) \oplus \sll(11)$, $\sll(2) \oplus \sll(41) \oplus \sll(11)$, etc. (cf. Remark 1 in \S \ref{sec:classification}).

To prove Theorem \ref{main thm}, we establish a one-to-one correspondence 
between complex IFPSs and prehomogeneous vector spaces, which is a purely algebraic 
object introduced by M.Sato \cite{sato}. 
He called a triplet $(G, \rho, V)$ a prehomogeneous vector space (abbrev. PV) if $G$ is a connected linear algebraic group over 
algebraically closed field $K$ and $\rho$ is a rational 
representation of $G$ on a finite dimensional $K$-vector space $V$ such that $V$ admits a Zariski-open $G$-orbit.  
When $K$ is equal to the complex number field, we can naturally extend the notion of PV to the holomorphic category in our sense (see \S 5). 

For a PV $(G, \rho, V)$, its infinitesimal form $(\g, d\rho, V)$ satisfies the condition that there exists 
$v \in V$ such that $d\rho(\g)v = V$ (cf. \cite{kimura}). We can show that a complex IFPS induces a Lie algebra representation satisfying the 
same condition as that of PVs. In this context we prove that IFPSs correspond to PVs.  

The paper is organized as follows. The first half of the paper, which is composed of $\S 2$-$\S 4$, is a geometric preliminary for 
proving the correspondence between IFPSs and PVs. In $\S 5$ we introduce the notion of PV, and prove this correspondence.  
In $\S 6$ we prepare some important notions of PV concerning classifications of PVs. 
Finally in $\S 7$ we prove Theorem \ref{main thm} by using a classification of certain irreducible PVs by Sato and Kimura \cite{sato-kimura}. 
\section{$(G, X)$-structures and flat Cartan structures}
A flat projective structure on a complex manifold $M$ is a special case of ($G$, $X$)-structures, which will be defined in the following.  
Let $G$ be a complex Lie group, and let $X$ be a connected complex homogeneous space of $G$.     
Then by the Liouville theorem (cf. [3, Proposition 1.5.2]), the action of $G$ on $X$ is locally determined. Namely, for $g \in G$ if there exists a nonempty open subset $U$ of $X$ such that $g$ gives the identity transformation of $U$, then $g$ gives the global identity transformation of $X$. 
From now on we assume that $G$ acts on $X$ effectively. Thus if $g$ gives a local identity transformation of $X$,  then $g$ is equal to the unit element of $G$. 
Accordingly we identify an element $g$ of $G$ with the transformation on some open set of $X$ induced by $g$.    
We assume that $\dim M = \dim X$.  
\renewcommand{\theenumi}{\arabic{enumi}}
\renewcommand{\labelenumi}{\rm (\theenumi)}
\begin{df}[\cite{goldman}]
A $(G, X)$-structure on $M$ is a maximal atlas $\{($$U_\alpha$, $\varphi_\alpha$$)\}_{\alpha \in A}$ of $M$ such that
\begin{enumerate}
\item $\varphi_\alpha$ maps $U_\alpha$ biholomorphically onto an open subset of $X$. 
\item For every pair $(\alpha, \beta)$ with $U_\alpha \cap U_\beta \neq \emptyset$ and each connected component $C$ of $U_\alpha \cap U_\beta$, 
there exists $\tau(C; \beta, \alpha) \in G$ such that 
$\varphi_\beta \circ \varphi_\alpha^{-1}|_{\varphi_\alpha(C)}$ equals the restriction of $\tau(C; \beta, \alpha)$ to $\varphi_\alpha (C)$.   
\end{enumerate}
\end{df}
Two atlases are said to be equivalent if they are compatible. Then note that 
a maximal atlas corresponds to an equivalence class of an atlas.  

In the following we fix a point $x \in X$, and denote the isotropy subgroup at $x$ by $G'$. Then we can identify $X$ with the quotient space
$G/G'$. 
In order to prove the correspondence between IFPSs and PVs, in this section we generally show that $(G, X)$-structures on $M$ correspond to 
flat Cartan structures of type $G/G'$ on $M$, which will be defined as follows:  
We denote the Lie algebra of $G$ by $\g$, and the Lie algebra of $G'$ by $\g'$.
\begin{df}[\cite{kobayashi}]
Let $\pi_P : P \to M$ be a principal $G'$ bundle,  and let $\omega$ be a $\g$-valued 1-form on $P$. 
We say that $(P, \omega)$ is a Cartan structure of type $G/G'$ on $M$ if
\begin{enumerate}
\item $R_a^* \omega = \Adj(a^{-1}) \omega$ \quad for $a \in G'$,
\item $\omega (A^*) =A$ \quad for $A \in \g'$, \\
where $A^*$ is the fundamental vector field corresponding to $A$.
\item For $u \in P$, $\omega_u : T_uP \to \g$ gives a linear isomorphism.
\end{enumerate}
\end{df}
A 1-form $\omega$ of a Cartan structure $(P, \omega)$ is called a Cartan connection. A Cartan structure $(P, \omega)$ is said to be flat 
if the equality $d\omega + \frac{1}{2}[\omega, \omega]=0$ holds.
Cartan structures $(P, \omega)$ and $(P', \omega')$ on $M$ are said to be isomorphic  (via identity transformation of $M$) if there exists a bundle isomorphism $\phi: P \to P'$ such 
that $\phi^* \omega' = \omega$, and $\phi$ induces the identity transformation of $M$. 
We call $\phi$ an isomorphism of Cartan structures. 
We denote this equivalence relation by $(P, \omega)$ $\sim$ $(P', \omega')$. 

Now $(G, X)$-structures correspond to isomorphism classes of flat Cartan structures $(P, \omega)$ as follows:
\begin{thm}\label{bij1}
There is a one-to-one correspondence between the set of $(G, X)$-structures on $M$ and the set of isomorphism classes of flat Cartan 
structures of type $G/G'$ on $M$. 
\end{thm}
\pf 
We shall construct a map
\begin{eqnarray*}
&&\Phi: \{\mbox{$(G, X)$-structure on $M$}\} \to  \hspace{3cm}\\
&& \quad \{\mbox{flat Cartan structure of type $G/G'$ on $M$}\}/\sim. 
\end{eqnarray*}
Let $\{(U_\alpha, \varphi_\alpha)\}_{\alpha \in A}$ be a $(G, X)$-structure on $M$.
First we construct a principal fiber bundle $P$.
As a chart of $M$ let us choose $\varphi_\alpha: U_\alpha \to \varphi_\alpha(U_\alpha) \subset X$. Then 
the principal fiber bundle $\pi: G \to X= G/G'$ induces a principal fiber bundle $\pi: \pi^{-1}(\varphi_\alpha(U_\alpha)) \to 
\varphi_\alpha(U_\alpha)$. We denote $\pi^{-1}(\varphi_\alpha(U_\alpha))$ by $\widetilde{U}_\alpha$. Through $\varphi_\alpha$ we obtain a 
principal fiber bundle $\pi_\alpha: \widetilde{U}_\alpha \to U_\alpha$. 
Let $\omega$ be the Maurer-Cartan form of $G$. Then $\omega$ induces a Cartan connection $\omega_\alpha$ on $\widetilde{U}_\alpha$ 
by setting $\omega_\alpha := \omega|_{\widetilde{U}_\alpha}$.
Hence we obtain a family of Cartan structures $\{(\widetilde{U}_\alpha, \omega_\alpha)\}_{\alpha \in A}$.

We denote by $\tau(C; \beta, \alpha)$ a coordinate change $\varphi_\beta \circ \varphi_\alpha^{-1}$ on a connected component $C$ 
of $U_\alpha \cap U_\beta$.  
Let $U_\gamma$ be another coordinate neighborhood, and let $D$ and $E$ be connected components of 
$U_\beta \cap U_\gamma$ and $U_\gamma \cap U_\alpha$ respectively. If $C \cap D \cap E \neq \emptyset$, then 
we have $\tau(D; \gamma, \beta) \tau(C; \beta, \alpha) = \tau(E; \gamma, \alpha)$. This follows from the fact that 
the action of $G$ on $X$ is locally determined (see the beginning of this section).    
Let $g \in \widetilde{U}_\alpha$, and $h \in \widetilde{U}_\beta$. We express $g \sim h$ if $\pi_\alpha 
(g) = \pi_\beta (h)$ and $h = \tau(C; \beta, \alpha) g$ on a connected component $C$ of $U_\alpha \cap U_\beta$ containing $\pi_\alpha (g)$. 
Then this is an equivalence relation in the set 
$\bigsqcup_{\alpha \in A} \widetilde{U}_\alpha$. 
Hence we obtain a quotient space $P := \bigsqcup_{\alpha \in A} \widetilde{U}_\alpha/\sim$. 
Then note that we have the natural inclusion $\imath$: $\widetilde{U}_\alpha \to P$.     

We show that $P$ is a principal $G'$ bundle over $M$. Let us define a projection $\pi_P: P \to M$ by $\pi_P(\imath(g)) = \pi_\alpha(g)$ for 
$g \in \widetilde{U}_\alpha$, and define a 
group 
action of $G'$ on $P$ by $\imath(g) a := \imath(ga)$ for $a \in G'$. A local trivialization of $P$ is derived from that of $\widetilde{U}_\alpha$ with its 
complex structure. These are well defined and satisfy the conditions of principal fiber bundles. 

Secondly we construct a flat Cartan connection $\omega_P$ on $P$. 
Here note that any tangent vector at $\imath(g)$ of $P$ is given by $d\imath (X)$ for some $X \in T_g\widetilde{U}_\alpha$. 
Now we define a 1-form $\omega_P$ on $P$ by $\omega_P (d\imath(X)) := \omega_\alpha(X)$ for $X \in T_g\widetilde{U}_\alpha$. Then $\omega_P$ 
is well defined. Since $\omega_\alpha$ is a flat Cartan connection on $\widetilde{U}_\alpha$, $\omega_P$ gives a flat Cartan connection on $P$.
As a result we have obtained a flat Cartan structure  $(P, \omega_P)$ of type $G/G'$ on $M$.

Finally we show that equivalent atlases of $(G, X)$-structure induce isomorphic Cartan structures.
Let two atlases $\{(U_\alpha, \varphi_\alpha)\}_{\alpha \in A}$ and $\{(V_\lambda, \psi_\lambda)\}_{\lambda \in \Lambda}$ belong to 
the same 
$(G, X)$-structure. Then each atlas on $M$ induces a Cartan structures $(P, \omega_P)$ and $(P', \omega_{P'})$ respectively.
We define a bundle map $f: P \to P'$ by $f(\imath(g)) = \imath(\tau(C; \lambda, \alpha) \cdot g)$ for $g \in \widetilde{U}_\alpha$ and 
$\pi_\alpha(g) \in V_\lambda$ such that $\pi_\alpha(g) \in C \subset U_\alpha \cap V_\lambda$. 
Then $f$ is well defined, and moreover $f$ is a holomorphic bundle isomorphism. We can easily verify that 
$f^* \omega_P' = \omega_P$, and that $f$ induces the identity transformation of $M$. 
Hence $(P, \omega_P)$ is isomorphic to $(P', \omega_{P'})$, and consequently we obtain a map 
$\Phi$: $\{$$(G, X)$-structure on $M$$\}$ $\to$ $\{$flat Cartan structure of type $G/G'$ on $M$$\}/\sim$.
\\[-2mm]

Next we shall construct a map
\begin{eqnarray*}
&&\Psi: \{\mbox{flat Cartan structure of type $G/G'$ on $M$}\}/\sim \ \to \\  
&& \quad \{\mbox{($G$, $X$)-structure on $M$}\}. 
\end{eqnarray*}
Let $(P, \omega_P)$ be a flat Cartan structure of type $G/G'$ on $M$. Let $\pi_P$: $P \to M$ be its projection. 
The Maurer-Cartan form $\omega$ of $G$ is a $\g$-valued  
1-form satisfying the structure equation $d\omega +\frac{1}{2}[\omega, \omega]=0$, which $\omega_P$ also satisfies. 
Hence for any $u \in P$, there exists a neighborhood $U$ of $\pi_P(u)$ and a bundle isomorphism $\tilde{f}:$ ${\pi_P}^{-1}(U) \to V$, 
where $V$ is an open subset of $G$,  
such that $\tilde{f}^* \omega = \omega_P$ (cf. \cite{cap-slovak}, p.74).  Therefore $\widetilde{f}$ induces 
a biholomorphic mapping of base 
spaces $f:$ $U \to \pi(V)$. 
Suppose that $\tilde{f'}$ is another bundle isomorphism $\tilde{f'}:$ ${\pi_P}^{-1}(U') \to V'$ such that $\pi_P(u) \in 
U'$, and $\tilde{f'}^* \omega = \omega_P$.  Let $f': U' \to \pi(V')$ be its inducing map.  
For an element $g$ of $G$, We denote its left translation by $L_g$.  
Then for a connected component $C$ of $U \cap U'$, 
there exists a unique element $g \in G$ such that $\tilde{f}' = L_g \circ \tilde{f}$ on   ${\pi_P}^{-1}(C)$ (cf. \cite{cap-slovak}, p.74). 
Hence by setting $\tau(C; U', U) := g$, the equality $f' \circ f^{-1} = \tau(C; U', U)$ holds on $f(C)$. Consequently we obtain an atlas $\{(U, f)\}$ of $(G, X)$-structure on $M$.  

Next we suppose that flat Cartan structures $(P, \omega_{P})$ and $(P', \omega_{P'})$ are isomorphic. Then we can verify that two atlases 
of $(G, X)$-structure on $M$ induced by $(P, \omega_{P})$ and $(P', \omega_{P'})$ are equivalent. Hence we obtain a map $\Psi$:  
$\{$flat Cartan structure of type $G/G'$ on $M$$\}/\sim$ $\to$ $\{$$(G, X)$-structure on $M$$\}$.\\[-2mm]

Two maps $\Phi$ and $\Psi$ obtained above satisfy $\Phi \circ \Psi = \Psi \circ \Phi = id$. Hence Theorem \ref{bij1} has been proved.
\hfill $\Box$ \\

{\it Remark.}  \
Let $M$ be a complex manifold, and let $\psi$ be a biholomorphic map of $M$. 
Let $(P, \omega)$ and $(P', \omega')$ be flat Cartan structures of type $G/G'$ on $M$, and let 
$\{(U_\alpha, \varphi_\alpha)\}_{\alpha \in A}$ and $\{(U'_\lambda, \varphi'_\lambda)\}_{\lambda \in \Lambda}$ be the induced 
($G$, $X$)-structures on $M$ respectively. We say that $(P, \omega)$ and $(P', \omega')$ are isomorphic via $\psi$ if there exists a bundle 
isomorphism 
$\tilde{\psi}: P \to P'$ such that $\tilde{\psi}^* \omega' =\omega$ which induces $\psi$. On the other hand we say that $\{(U_\alpha, \varphi_\alpha)\}_{\alpha \in A}$ and 
$\{(U'_\lambda, \varphi'_\lambda)\}_{\lambda \in \Lambda}$ are equivalent via $\psi$ if $\psi: M \to M$  
satisfies the following condition: 
for each connected component $C$ of $U_\alpha \cap \psi^{-1}(U'_\lambda)$ there exists an element $g$ of $G$ such that the composite 
$\varphi'_\lambda \circ \psi \circ \varphi_\alpha^{-1}$ equals $g$ on $\varphi_\alpha(C)$.  
Then we can verify that $(P, \omega)$ and $(P', \omega')$ are isomorphic via $\psi$ if and only if 
$\{(U_\alpha, \varphi_\alpha)\}_{\alpha \in A}$ and $\{(U'_\lambda, \varphi'_\lambda)\}_{\lambda \in \Lambda}$ are equivalent via $\psi$. 
\\[-2mm]

Next we shall consider left invariant $(G, X)$-structures on a complex Lie group $L$. For an element $a$ of $L$ we denote the 
left translation of $L$ by $L_a$.   
\begin{df}
A $(G, X)$-structure $\{(U_\alpha, \varphi_\alpha)\}_{\alpha \in A}$ on $L$ is said to be left invariant if 
for any $a \in L$ and each connected component $C$ of 
$U_\alpha \cap L_a^{-1}(U_\beta)$, there exists an element $g$ of $G$ such that the composite $\varphi_\beta \circ L_a \circ \varphi_\alpha^{-1}|_{\varphi_\alpha (C)}$ 
equals the restriction of $g$ to $\varphi_\alpha (C)$.    
\end{df}

In terms of Cartan structures, left invariance is described by the following: 
\begin{df}[\cite{mendez-lopera}] 
Let $(P, \omega)$ be a Cartan structure of type $G/G'$ on $L$.  
Then $(P, \omega)$ is said to be left invariant if there exists a left action $L \times P \ni (a, u) \to {L_a}'(u) \in P$ 
satisfying the following conditions: For any $a \in L$, ${L_a}'$ is a bundle isomorphism such that  
$\pi_P \circ {L_a}' = L_a \circ \pi_P$ and ${L_a'} ^* \omega = \omega$. 
\end{df}
Let $(P, \omega)$ and $(P', \omega')$ be left invariant Cartan structures on $L$. We suppose that they are isomorphic as Cartan structures, thus 
there exists an isomorphism of Cartan structures $\phi$: $(P, \omega) \to (P', \omega')$.  
Then note that we have $\phi \circ {L_a'} = {L_a'} \circ \phi$ for $a \in L$.  
\begin{thm}\label{bij2}
There is a one-to-one correspondence between the set of left invariant $(G, X)$-structures on $L$ and the set of isomorphism classes 
of left invariant flat Cartan structures of type $G/G'$ on $L$. 
\end{thm}
\pf 
We show that the maps $\Phi$ and $\Psi$ in the proof of Theorem \ref{bij1} preserve left invariance. 

Let $\{(U_\alpha, \varphi_\alpha)\}_{\alpha \in A}$ be an atlas of  left invariant $(G, X)$-structure on $L$. 
We suppose that in Theorem \ref{bij1}, $\Phi(\{(U_\alpha, \varphi_\alpha)\})$ is given by an isomorphism class of $(P, \omega_P)$, 
where  $P = \bigsqcup_{\alpha \in A} \widetilde{U}_\alpha/\sim$ and $\omega_P (d\imath(X)) = \omega_\alpha(X)$ 
for $X \in T_g\widetilde{U}_\alpha$.
Then for $a \in L$, the left translation $L_a$ induces a bundle isomorphism ${L_a}'$: $P \to P$ defined by ${L_a}'(\imath(g)) = 
\imath((\varphi_\beta \circ L_a \circ {\varphi_\alpha}^{-1}) \cdot g)$ for $g \in \widetilde{U}_\alpha$ and $\beta \in A$ such that $a 
\pi_\alpha(g) \in 
U_\beta$. Indeed this definition is well defined. Furthermore we can verify that 
$\pi_P \circ {L_a'} = L_a \circ \pi_P$, and ${L_a'} ^* \omega = \omega$. 
We can also show that equivalent atlases of left invariant $(G, X)$-structure on $L$ induce isomorphic left invariant Cartan structures on 
$L$. Consequently $\Phi$ preserves left invariance.\\[-2mm]

Next we show that $\Psi$ preserves left invariance.
Let $(P, \omega_P)$ be a left invariant flat Cartan structure of type $G/G'$ on $L$. 
Recall that in the proof of Theorem \ref{bij1} 
$(P, \omega_P)$ induces a set $\{\tilde{f_\alpha}: {\pi_P}^{-1}(U_\alpha) \to G\}_{\alpha \in A}$ such that 
$\bigcup_\alpha U_{\alpha \in A} = L$. Then this set induces an atlas $\{(U_\alpha, f_\alpha)\}_{\alpha \in A}$ of $(G, X)$-structure on 
$L$.  For any $a \in L$ there is a bundle isomorphism ${L_a}': P \to P$ such that ${L_a'}^* \omega_P = \omega_P$, therefore  
$\tilde{f_\beta} \circ {L_a'} \circ \tilde{f_\alpha}^{-1}$ preserves the Maurer-Cartan form $\omega$ of $G$.
Hence for each connected component $C$ of $U_\alpha \cap L_a^{-1}(U_\beta)$, the left translation $\tilde{f_\beta} \circ {L_a'} \circ \tilde{f_\alpha}^{-1}$ is  given by a unique element $g$ of $G$ on $\tilde{f_\alpha}({\pi_P}^{-1}(C))$. From this fact, we can verify that 
the induced map $f_\beta \circ L_a \circ f_\alpha^{-1}$ is given by $g$ on $f_\alpha(C)$. 

Therefore the atlas $\{(U_\alpha, f_\alpha)\}_{\alpha \in A}$ is left invariant.
Since left invariance is preserved in the equivalence relation of atlases, we can conclude that $\Psi$ preserves left invariance. \\[-2mm]

We have proved that two maps $\Phi$ and $\Psi$ preserve left invariance. Since $\Psi$ is an inverse map of $\Phi$, Theorem \ref{bij2} 
has been proved.   
\hfill $\Box$ \\

\remark
Theorems \ref{bij1} and \ref{bij2} are valid also in the real $C^\infty$ category.  
Let $N$ be a real manifold. Then we can consider ($G$, $X$)-structures and Cartan structures similarly to the complex case.  
Concerning Theorem \ref{bij1}, in \cite[p.75]{cap-slovak} it has been stated that a ($G$, $X$)-structure exists on $N$ if and 
only if there exists a flat Cartan structure of type $G/G'$ on $N$.  

Here let us consider the projective geometry ($PGL(\R^{n+1})$, $P(\R^{n+1})$) for example. Then $P(\R^{n+1})$ is connected and 
$PGL(\R^{n+1})$ acts on $P(\R^{n+1})$ effectively. 
Thus we can apply these theorems to ($PGL(\R^{n+1})$, $P(\R^{n+1})$)-structures on a manifold $N$.  
Let us fix the point $o=[0, \cdots, 0, 1]$ $\in$ $P(\R^{n+1})$, and denote the isotropy subgroup at $o$ by $PGL(\R^{n+1})_o$.  
Then the real projective space $P(\R^{n+1})$ is expressed as $PGL(\R^{n+1})/PGL(\R^{n+1})_o$.  
 
In our definition a flat projective structure on $N$ is a ($PGL(\R^{n+1})$, $P(\R^{n+1})$)-structure on $N$, 
which corresponds to an 
isomorphism class of a flat Cartan structure $(P, \omega)$ of type $PGL(\R^{n+1})/PGL(\R^{n+1})_o$ on $N$. 

On the other hand there is a more familiar definition of flat projective structures, which is defined to be a projective equivalence class of
a projectively flat torsionfree affine connection.  In this case we can see that there is a one-to-one correspondence between the set of 
projective equivalence classes of projectively flat torsionfree affine connections 
on $N$ and the set of isomorphism classes of flat Cartan structures of type $PGL(\R^{n+1})/PGL(\R^{n+1})_o$ on $N$ by using the results of Tanaka \cite{tanaka} which is quoted 
in \cite[p.131]{agaoka}. Hence there is also a one-to-one correspondence between the set of ($PGL(\R^{n+1})$, $P(\R^{n+1})$)-structures on $N$ and 
the set of projective equivalence classes of projectively flat torsionfree affine connections on $N$. 
That is the reason why we call a ($PGL(\R^{n+1})$, $P(\R^{n+1})$)-structure a flat projective structure. 

\section{Left invariant flat Cartan structures and transitive embeddings}
Let $L$ be a complex Lie group, and let $\l$ be its Lie algebra.  
In this section we prove that left invariant flat Cartan structures of type $G/G'$ on $L$ correspond to certain injective Lie algebra homomorphisms of $\l$ to $\g$, which we call simply transitive embeddings (cf. Theorem 2.3 in \cite{mendez-lopera}). 
Recall that $\dim \l = \dim \g/\g'$. 
\begin{df}
Let $f : \l \to \g$ be a complex Lie algebra homomorphism. Then we call $f$ a simply transitive embedding {\rm(}of $\l$ into $\g${\rm)} if the induced linear map 
$\bar{f}: \l \to \g/\g'$ is a linear isomorphism. 
\end{df}
\remark
Let $\g_{-1}$ be a complementary subspace of $\g'$ in $\g$. We denote by $f_{-1}$ the $\g_{-1}$-component of $f$.
Let $f : \l \to \g$ be a complex Lie algebra homomorphism. Then $f$ is a simply transitive embedding if and only if 
$f_{-1} (\l) = \g_{-1}$. \\[-2mm]

Let $L$ be a real Lie group, and let $\l$ be its Lie algebra. Then  
Mendez, Lopera (Theorem 2.2 in \cite{mendez-lopera}) showed that a left invariant flat real Cartan structure 
$(P, \omega)$ of type 
$G/G'$ on $L$ induces a real Lie algebra homomorphism $f: \l \to \g$. The construction of 
homomorphism $f$ is essentially taken from \cite{agaoka}.   
In the holomorphic category, we have the same assertion.  In the following we summarize its proof: 
Let $L$ be a complex Lie group, and let $\l$ be its Lie algebra. 
Let $(P, \omega)$ be a flat (complex) Cartan structure of type $G/G'$ on a complex Lie group $L$. 
Since $(P, \omega)$ is left invariant, there is a left action of $L$ on $P$, i.e. $L \times P$ $\ni$ $(a, u) \mapsto L_a'(u)$ $\in P$. 
Let $\pi_P: P \to L$ be the projection of $(P, \omega)$. 
We fix an element $\hat{o}$ of $P$ such that $\pi_P(\hat{o})$ is equal to the unit element $e$ of $L$. 
We define a map $j: L \to P$ by $j(a) = L_a'(\hat{o})$.  Now let us consider 
the $\g$-valued 1-form $j^* \omega$ on $L$.  The 1-form $j^* \omega$ is left invariant,   
since $j$ is compatible with the left action of $L$ and $L_a'^* \omega = \omega$. Thus we obtain a linear map 
$j^* \omega: \l \to \g$. We denote this linear map by $f$. Then $f$ is a Lie algebra homomorphism if and 
only if $(P, \omega)$ is flat. 

Furthermore we can show that $f$ is a simply transitive embedding.  
\begin{prop}\label{cartan and (G)-hom}
Any left invariant flat Cartan structure $(P, \omega)$ of type $G/G'$ on $L$ induces a simply 
transitive embedding $f: \l \to \g$.
\end{prop}
\pf 
As we have seen above  $(P, \omega)$ induces a Lie algebra homomorphism $f=j^* \omega: \l \to \g$ by Theorem 2.2 
in \cite{mendez-lopera}. 
Since the map $j: L \to P$ satisfies $\pi \circ j = id$,  we have $j_*(T_eL) \oplus T_{\hat{o}}\pi^{-1}(\hat{o}) = T_{\hat{o}}P$. Since the Cartan connection 
$\omega$ gives 
a linear isomorphism $\omega_{\hat{o}}: T_{\hat{o}}P \to \g$, we have $\omega (j_*(T_eL)) \oplus \g' = \g$. Hence $f(\l$) is isomorphic to 
$\g/\g'$, which implies that the induced linear map $\bar{f}: \l \to \g/\g'$ is an isomorphism. 
Consequently $f$ is a simply transitive embedding.  \hfill $\Box$

\begin{df}\label{equivalence of (G)-hom}
We say that simply transitive embeddings $f_1$ and $f_2$ are equivalent 
if there exists $g \in G'$ such that $f_2 = \Adj(g) f_1$. We denote this equivalence relation by $f_1 \sim f_2$.  
\end{df}
We denote the equivalence class of a simply transitive embedding $f$ by $[f]$.  
\begin{lem}\label{(G)-hom}
Isomorphic left invariant flat Cartan structures of type $G/G'$ on $L$ induce equivalent simply 
transitive embeddings. 
\end{lem}
\pf 
Let $(P_i, \omega_i)$ be a left invariant flat Cartan structure of type $G/G'$ on $L$  
($i = 1, 2$).  Then by Proposition \ref{cartan and (G)-hom} each $(P_i, \omega_i)$ induces a simply transitive embedding $f_i$. 
Now we assume that $(P_1, \omega_1)$ is isomorphic to $(P_2, \omega_2)$. Then there exists a bundle isomorphism $\phi: P_1 \to P_2$ 
such that 
$\phi^* \omega_2 = \omega_1$, and $\phi \circ L_a' = L_a' \circ \phi$. Let $\hat{o}_i$ be a fixed element of $P_i$ such that $\pi (\hat{o}_i) =e$. Then 
there exists $g \in G'$ such that $\phi(\hat{o}_1) = \hat{o}_2 \cdot g$. 
The left action of $L$ on $P_i$ induces a map $j_i: L \to P_i$ by $j_i(a) = L_a' (\hat{o}_i)$. Then we have $\phi \circ j_1 = R_g \circ j_2$. 
Hence we have 
$f_1 = j_1^* \omega_1 = j_1^* \phi^* \omega_2 = (R_g \circ j_2)^* \omega_2 = \Adj(g^{-1}) j_2^*\omega_2 =\Adj(g^{-1})f_2$. \hfill $\Box$ \\[-3mm]

By combining Proposition \ref{cartan and (G)-hom} and Lemma \ref{(G)-hom} we obtain a map
\begin{eqnarray*}
&&\Theta: \{\mbox{left invariant flat Cartan structure of type $G/G'$ on $L$}\}/\sim \\ 
&& \quad \to \{\mbox{simply transitive embedding} \  f: \l \to \g\}/\sim. 
\end{eqnarray*}
We denote the equivalence class of a Cartan structure $(P, \omega)$ by $[(P, \omega)]$.   
\begin{prop}\label{bij3}
$\Theta$ is bijective.
\end{prop}

\pf
First we show that $\Theta$ is injective. 
Let $(P_i, \omega_i)$ be a left invariant flat Cartan structure of type $G/G'$ on $L$ ($i = 1, 2$). 
We suppose that $\Theta([(P_1, \omega_1)]) = \Theta([(P_2, \omega_2)])$. Then by the definition of $\Theta$ there exists $g \in G'$ such that 
$j_2^* \omega_2 = \Adj(g) j_1^* \omega_1$. 
We define a map $\phi: P_1 \to P_2$ by $\phi(j_1(a)h) = j_2(a)gh$ for $h \in G'$. We prove that $\phi$ gives an isomorphism between 
$(P_1, \omega_1)$ and $(P_2, \omega_2)$. 
We can easily verify that $\phi \circ j_1 = R_g \circ j_2$, and $\phi \circ 
L_a' = L_a' \circ \phi$ for $a \in L$. Next we show that $\phi$ preserves Cartan connections. Any tangent vector of $P_1$ can be uniquely 
expressed in the form ${R_h}_*({j_1}_*X + Z^*)$ where $X \in T_aL$, $Z \in \g'$ and $h \in G'$. Then 
$\phi^* \omega_2 ({j_1}_* X) =\omega_2({R_g}_* {j_2}_* X) = \Adj(g^{-1}) j_2^* \omega_2 (X) = \Adj(g^{-1}) \Adj(g) j_1^*\omega_1(X) = 
\omega_1({j_1}_*X)$. 
Moreover for $u \in P$, $\phi^* \omega_2 (Z^*_u) = \omega_2 (\phi_* Z^*_u) = \omega_2 (Z^*_{\phi(u)}) 
= Z = \omega_1 (Z^*_u)$. Hence \\[-6mm]
\begin{eqnarray*}
\phi^* \omega_2 ({R_h}_*({j_1}_*X + Z^*)) &=& \Adj(h^{-1}) \phi^* \omega_2 ({j_1}_*X +Z^*) \\ 
&=& \Adj(h^{-1}) \omega_1 ({j_1}_*X +Z^*) \\
&=& \omega_1 ({R_h}_* (j_*X +Z^*)). 
\end{eqnarray*}
It follows that $\phi^* \omega_2 = \omega_1$, which implies 
that $(P_1 ,\omega_1)$ is isomorphic to $(P_2, \omega_2)$. \\[-1mm]

Secondly we show that $\Theta$ is surjective. 
Let $f: \l \to \g$ be a simply transitive embedding.  
We shall construct a map $j: L \to P$. 
Fix a linear frame $\tilde{o}$ at the unit element $e$ of $L$. Then a left invariant frame field $\{{L_a}_* \tilde{o}\}_{a \in L}$ gives a 
complete 
parallelism on $L$. Hence $L$ has an $\{e\}$-structure $\tilde{L}$, i.e. an $\{e\}$-reduction of the frame bundle of $L$. 
Let $h: \tilde{L} \to P$ be 
an extension of $\tilde{L}$ by the injective homomorphism $\{e\} \to G'$. We define a map $\tilde{j}: L \to \tilde{L}$ by $\tilde{j}(a) = {L_a}_* 
(\tilde{o})$. 
We denote a composite $h \circ \tilde{j}$ by $j$. Then for the natural projection $\pi$ we have $\pi \circ j = id$. 
By using the map $j$, we construct a left invariant flat Cartan connection $\omega$ on $P$ such that $j^* \omega =f$, 
following the proof of Theorem 2.12 in \cite{agaoka}. 
For $a \in L$ any tangent vector at $j(a)$ can be uniquely written in the form $j_*X + A^*$ where $X \in T_aL$, and $A \in \g'$. 
We set $\omega_{j(a)} (j_*X +A^*):= f(X) +A$, and extend it to any point of $P$ by $\omega_{j(a)g} := \Adj(g^{-1}) R_{g^{-1}}^* 
\omega_{j(a)}$ 
for $a \in L$ and $g \in G'$. Then $(P, \omega)$ gives a flat Cartan structure of type $G/G'$, and obviously $j^* \omega =f$. 

Next we show that $(P, \omega)$ is left invariant. Any point $u \in P$ is uniquely expressed in the form $u=j(b)g$ 
$(b \in L, g \in G')$. We define a map $L_a': P \to P$ by $L_a'(j(b)g) = j(ab)g$. Then $L_a'$ gives a bundle isomorphism of $P$. ${L_a'}$ 
defines a left action of $L$ on $P$ and satisfies $\pi \circ L_a' = L_a \circ \pi$. 
We can easily verify that $L_a^* \omega = \omega$ for $a \in L$. Hence $(P, \omega)$ is left invariant. 
If we set $h(\tilde{o}) = \hat{o}$, then $(P, \omega)$ induces $f$. Therefore $\Theta$ is surjective. 
\hfill $\Box$
\\[-1mm]

\remark
Mendez, Lopera \cite[Theorem 2.3]{mendez-lopera} proved the following result: Let 
$G$ be a connected and simply connected real Lie group, and let $H$ be a connected closed subgroup. Let $A/B$ a real homogeneous 
space satisfying $\dim G/H = \dim A/B$. Let $\g$ be the Lie algebra of $G$, and let $\da$ be the Lie algebra of $A$, and let $\db$ be the Lie algebra of $B$. 
Then 
they proved that there exists an invariant flat Cartan structure of type $A/B$ on the homogeneous 
space $G/H$ 
if and only if there exists a Lie algebra homomorphism $f: \g \to \da$ such that 
$f(\h) \subset \db$ and the induced map $\hat{f}: \g/\h \to \da/\db$ is an isomorphism.\\[-3mm]

By combining Theorem \ref{bij2} and Proposition \ref{bij3} we obtain the following Theorem: 
\begin{thm}\label{(G, X)-str and (G)-hom} 
There is a one-to-one correspondence between the set of left invariant $(G, X)$-structures on a complex Lie group $L$ and the set of 
equivalence classes of simply transitive embeddings of $\l$ into $\g$.  \\[-4mm]
\end{thm}

{\it Remark 1.}  \ 
Kim \cite[Theorem 2.4]{kim} proved a similar result: A connected and simply connected real Lie group $L$ admits a $(G, X)$-structure if and only if 
there exists a Lie algebra 
homomorphism $f: \l \to \g$ satisfying the following condition: 
for the isotropy subalgebra $\g_x$ at some $x \in X$, we have $f(\l) \cap \g_x =0$.  
This condition is equivalent to the one that $f$ is a Lie algebra homomorphism of $\l$ 
such that $\bar{f}: \l \to \g/\g_x$ is a linear isomorphism for some $x \in X$. \\[-2mm] 

{\it Remark 2.} \  On a real manifold we can consider $(G, X)$-structures similarly to the complex case. 
Let $G$ be a real Lie group, and $X$ be its connected homogeneous space. As in the complex case 
we suppose that $G$ acts on $X$ effectively.   
Then all the assertions we have proved are true in the real $C^\infty$ category. 
Here we explain the relationship between complex geometric structures and real geometric structures.
  
Let $G$ be a complex Lie group, and let $G'$ be its closed complex subgroup. We assume that $G/G'$ admits a real form 
$G_r/G_r'$. Namely, $G_r$ (resp. $G_r'$) is a real Lie group whose Lie algebra is a real form $\g_r$ (resp. $\g'_r$) of $\g$ (resp. $\g'$). 
Moreover we suppose that $G_r/G_r'$ is connected, and $G_r$ acts on  $G_r/G_r'$ effectively. 

Firstly let us consider a left invariant $(G, G/G')$-structure on a complex Lie group $L$. By Theorem \ref{(G, X)-str and (G)-hom} we have the 
corresponding simply transitive embedding $f: \l \to \g$. 
Assume that $\l$ has a real form $\l_r$, and let $L_r$ be a real Lie group whose Lie algebra is $\l_r$. 
Furthermore we suppose that $f$ has a real form
$f_r: \l_r \to \g_r$. Then $f_r$ gives a simply transitive embeding,   
and therefore $L_r$  admits a left invariant $(G_r, G_r/G'_r)$-structure. 

Conversely we suppose that a left invariant $(G_r, G_r/G'_r)$-structure is given on a real Lie group $L$. 
Let $f: \l \to \g_r$ be its corresponding simply transitive embedding.  
We denote by $\l^{\C}$ the complexification of $\l$. Let $L^{\C}$ be a complex Lie group having $\l^{\C}$ as its Lie algebra. 
Then the complexification 
$f^{\C}: \l^{\C} \to \g$ gives a simply transitive embedding. Hence $L^{\C}$ 
admits a left invariant $(G, G/G')$-structure. 

\section{IFPSs and corresponding Lie algebra homomorphisms} 
Let $L$ be a complex Lie group of dimension $n$. From now on we consider left invariant flat complex projective structures on $L$.  
We abbreviate a left invariant flat complex projective structure to complex IFPS. 
We assume that 
$G$ is the complex projective transformation group $PGL(\C^{n+1})=GL(\C^{n+1})/{\C^\times I}$, and $G'$ is the isotropy subgroup at the point 
$o$ $=$ [$0$, $\cdots$, $0$, $1$]  
of the  complex projective space $P(\C^{n+1})$. 
Note that $G/G'$ is connected and $G$ acts on $G/G'$ effectively. 
A $($$PGL(\C^{n+1})$, $P(\C^{n+1})$$)$-structure on a complex manifold is called a flat complex projective structure. 
We rewrite Theorem \ref{(G, X)-str and (G)-hom} for left invariant flat complex projective structures. 

Let $\g$ be the Lie algebra of $G$. Then $\g$ is isomorphic to $\sll(\C^{n+1})$, which is decomposed into $\g_{-1} \oplus \g'$ as follows:
\begin{eqnarray*}
&\g_{-1}&=
\left\{
      \left.
           \begin{pmatrix}
                 0 & u \\
                 0 & 0 
           \end{pmatrix}  
      \right|
          u \in \C^n
\right\},\\
&\g'&=
\left\{
     \left.
           \begin{pmatrix}
                 B & 0 \\
                 \xi & -\mathrm{tr}B 
           \end{pmatrix}  
      \right|
          B \in \mathfrak{gl}(\C^n), \ \xi \in {\C^n}^*
\right\}.
\end{eqnarray*}

\vspace{2mm}
\noindent
Let $f: \l \to \g$ be a Lie algebra homomorphism. Then 
$f$ is a simply transitive embedding if and only if $f$ satisfies   
$f(\l) e_{n+1} \oplus \langle e_{n+1} \rangle$ $=$ 
$\C^{n+1}$ where $e_{n+1}$ is the $(n+1)$-th vector of the standard basis. 
Then from Theorem \ref{(G, X)-str and (G)-hom} we obtain a one-to-one correspondence between the set 
$\{\mbox{complex IFPS on} L\}$ 
and the set   
$\{\mbox{Lie algebra homomorphism} f: \l \to \g \mid f(\l) e_{n+1} \oplus \langle e_{n+1} \rangle = \C^{n+1}\}/\sim.$ 
Agaoka proved almost the same result in the real case (\cite{agaoka}, Theorem 2.12). He defines a flat projective structure by using not atlases 
but 
linear connections. In \cite{agaoka} a Lie algebra homomorphism $f: \l \to \g$ corresponding to an IFPS is called a (P)-homomorphism. The 
condition of (P)-homomorphism is equivalent to the condition $f(\l) e_{n+1} \oplus \langle e_{n+1} \rangle$ $=$ 
$\R^{n+1}$.  
In the following we shall describe this condition in a little generalized form.  
We denote by $\pi$ the natural projection $\pi:$ $\gl(\C^{n+1}) \to \gl(\C^{n+1})/\C I_{n+1}$. For an element $P$ of $GL(\C^{n+1})$ 
(resp. $v$ of $\C^{n+1}$) we denote by $\overline{P}$ (resp. $\overline{v}$) its projection onto $PGL(\C^{n+1})$ (resp. $P(\C^{n+1})$). 

\begin{df}
Let $f, g$ be Lie algebra homomorphisms from $\l$ to $\gl(\C^{n+1})$, and let $v, w$ be vectors of $\C^{n+1}$. Then we say that 
$(f, v)$ and $(g, w)$ are equivalent if there exists $P \in GL(\C^{n+1})$ such that $\pi \circ g = \Adj(\overline{P}) (\pi \circ f)$ and 
$\overline{w} =\overline{P} \overline{v}$. We denote this equivalence relation by $(f, v) \sim (g, w)$. 
\end{df}

We denote the equivalence class of $(f, v)$ by $[(f, v)]$.  
In the following lemma an equivalence relation in the latter set 
is the one in Definition \ref{equivalence of (G)-hom}. 
\begin{lem}\label{lem2}
There is a one-to-one correspondence between the set 
\begin{eqnarray*}
&&\hspace{-0.3cm}\{(f, v) \mid f: \l \to \gl(\C^{n+1}) \ \mbox{is a Lie algebra homomorphism of $\l$}, v \in \C^{n+1} \\
&& \mbox{such that} \ f(\l)v \oplus \langle v \rangle = \C^{n+1}\}/\sim \\[2mm]
&& \hspace{-0.5cm} \mbox{and the set} \  \{f': \l \to \sll(\C^{n+1}) \mid f' \ \mbox{is a simply transitive embedding}\}/\sim. 
\end{eqnarray*}
\end{lem}
\pf 
Let $f: \l \to \gl(\C^{n+1})$ be a Lie algebra homomorphism and let $v$ be a vector of $\C^{n+1}$ such that 
$f(\l)v \oplus \langle v \rangle = \C^{n+1}$. There exists an element $P \in GL(\C^{n+1})$ such that $Pv=e_{n+1}$. Then    
the map $f' :=$ $\Adj(\overline{P})(\pi \circ f)$: $\l \to \sll(\C^{n+1})$ is a Lie algebra homomorphism and satisfies  
$f'(\l) e_{n+1} \oplus \langle e_{n+1} \rangle$ $=$ $\C^{n+1}$. Hence by setting $\sigma([f, v]) = [f']=[\Adj(\overline{P}) (\pi \circ f)]$, 
we obtain the map  
\begin{eqnarray*}
&&\hspace{-0.5cm}\sigma:\{(f, v) \mid \mbox{f is a Lie algebra homomorphism} \l \to \gl(\C^{n+1}), v \in \C^{n+1} \\
&& \quad \mbox{such that} \ f(\l)v \oplus \langle v \rangle = \C^{n+1}\}/\sim \quad \to \\ 
&& \{\mbox{homomorphism} \ f': \l \to \sll(\C^{n+1}) \mid f'(\l) e_{n+1} \oplus \langle e_{n+1} \rangle = \C^{n+1} \}/\sim.   
\end{eqnarray*}
This map $\sigma$ is well defined, and clearly surjective. We show that $\sigma$ is injective. 
Suppose that $\sigma([f, v]) = \sigma([g, w])$. By definition there exist $P$, $Q$ $\in GL(\C^{n+1})$ such that 
$\sigma([f, v]) = [\Adj(\overline{P}) (\pi \circ f)]$ and $\sigma([g, w]) = [\Adj(\overline{Q}) (\pi \circ g)]$. 
From the assumption we have $\overline{A} \in G'$ such that 
$\Adj(\overline{Q}) (\pi \circ g) = \Adj(\overline{A}) \Adj(\overline{P})(\pi \circ f)$. 
Hence we have $\pi \circ g = \Adj(\overline{Q^{-1} A P}) (\pi \circ f)$. 
We can easily check that $\overline{Q^{-1} A P} \overline{v} =\overline{w}$, and therefore $(g, w)$ is equivalent to $(f, v)$.
Hence $\sigma$ is injective. 
\hfill $\Box$
\\[-1mm]

From Theorem \ref{(G, X)-str and (G)-hom} and Lemma \ref{lem2}, we obtain the following: 
\begin{thm}\label{IFPS and hom}
There is a one-to-one correspondence between the set 
\begin{eqnarray*}
&&\hspace{-0.3cm}\{\mbox{complex IFPS on} \ L\} \ \mbox{and the set} \\ 
&&\hspace{-0.3cm}\{(f, v) \mid f: \l \to \gl(\C^{n+1}) \ \mbox{is a Lie algebra homomorphism and}\  v \in \C^{n+1} \\ 
&&\hspace{-0.3cm} \quad \mbox{such that} \ f(\l)v \oplus \langle v \rangle = \C^{n+1}\}/\sim. 
\end{eqnarray*}
\end{thm}
\remark
In the real case we have the same assertion. The corresponding  
condition  
$f(\l)v \oplus \langle v \rangle = \R^{n+1}$ has been obtained by Urakawa \cite[p.348]{urakawa}.  
\begin{df}\label{irreducibility}
Let $\{(U_\alpha, \phi_\alpha)\}_{\alpha \in A}$ be a complex IFPS on $L$. We suppose that $\{(U_\alpha, \phi_\alpha)\}_{\alpha 
\in A}$ corresponds to an equivalence class $[(f, v)]$ in Theorem {\rm \ref{IFPS and hom}}. 
We say that a complex IFPS $\{(U_\alpha, \phi_\alpha)\}_{\alpha \in A}$ is irreducible {\rm(}resp. reducible{\rm)}  if $f$ is irreducible 
{\rm(}resp. reducible{\rm)}. 
\end{df}

\section{Prehomogeneous vector spaces and IFPSs}\label{sec:PV and IFPS}
In this section, we explain that a complex IFPS on a Lie group corresponds to a certain prehomogeneous vector space.    
The notion of prehomogeneous vector space is originally introduced by Sato \cite{sato} in an algebraic category as follows: 
Let $G$ be a connected linear algebraic group over an algebraically closed field $K$, and let $\rho: G \to GL(V)$ be its finite dimensional 
rational representation on a $K$-vector space $V$. In this paper we call this triplet $(G, \rho, V)$ an algebraic triplet over $K$. 
In \cite{sato} the algebraic triplet $(G, \rho, V)$ is 
called a prehomogeneous vector space if $V$ has a Zariski-open $G$-orbit. 
We extend the notion of prehomogeneous vector spaces over the complex number field $\C$ to the holomorphic category as follows. 
In this paper we assume that an algebraic triplet always means an algebraic triplet over $\C$. 
We say that a subset  $O$ of $\C$-vector 
space $V$ is Euclidean-open if $O$ is open with respect to the Euclidean topology. Let $G$ be a complex Lie group, and let 
$\rho: G \to GL(V)$ be a finite dimensional holomorphic representation. We call this triplet $(G, \rho, V)$ a holomorphic triplet. 
\begin{df}
Let $(G, \rho, V)$ be a holomorphic triplet. 
We call the triplet $(G, \rho, V)$ a prehomogeneous vector space {\rm(}abbrev. PV{\rm)} if there exists $v \in V$ such that $\rho(G)v$
is Euclidean-open. 
\end{df}
Let $(G, \rho ,V)$ be an algebraic triplet over $\C$, and let $(\g, d\rho, V)$ be its infinitesimal form. 
Then $\rho(G)v$ is Zariski-open if and only if $d\rho(\g)v = V$ (cf. \cite{kimura}, Proposition 2.2). 
A point $v$ is said to be generic if $v$ belongs to a Zariski-open orbit. 
For a generic point $v \in V$, 
the isotropy subgroup $G_v$ at $v$ is called a 
generic isotropy subgroup. We denote by $\g_v$ its Lie algebra, i.e. $\g_v$ $=$ $Lie (G_v)$ $=$ $\{X \in \g \mid d\rho(X)v=0\}$. Then 
$d\rho(\g)v = V$ if and 
only if $\dim G - \dim G_v = \dim V$.   
These equivalences also hold in the holomorphic category, which we shall show in the following. 
\begin{prop}\label{prop1}
Let $(G, \rho, V)$ be a holomorphic triplet, and let $(\g, d\rho, V)$ be its infinitesimal form.
Then
\begin{enumerate}
\item $\rho(G)v$ is Euclidean-open if and only if $d\rho(\g)v = V$, 
\item $d\rho(\g)v = V$ if and only if $\dim \g - \dim \g_v = \dim V$. 
\end{enumerate}
\end{prop}
\pf 
The assertion (2) is obvious. We prove the assertion (1). 

We define a map $\phi: G \to V$ by $\phi(g) = \rho(g)v$. Then its differential $d\phi: \g \to V$ is given by $d\phi(X) = d\rho(X)v$.  
We suppose that $\rho(G)v$ is Euclidean-open in $V$. The map $\phi$ induces the map $\bar{\phi}: G/G_v \to V$ by 
$\bar{\phi}(gG_v) = \rho(g)v$. 
Its differential $d\bar{\phi}: \g/\g_v \to V$ is given by $d\bar{\phi}(X) = d\rho(X)v$.  
Since the homogeneous space $G/G_v$ is biholomorphic to $\rho(G)v$ through $\bar{\phi}$, the differential $d\bar{\phi}: \g/\g_v \to V$ is a 
linear isomorphism.  Hence we have $d\rho(\g)v =V$. 

Next we suppose that $d\rho(\g)v =V$. Then $d\phi_e: T_e G \to V$ is surjective. 
For any $g \in G$, the differential $d\phi_g: T_gG \to V$ is given by 
$d\phi_g =\rho(g) {L_g^{-1}}^* d\phi_e$, where ${L_g}^*$ is the pull back by the left translation $L_g$. 
Hence $d\phi_g$ is surjective, and there exists a neighborhood $U_g$ of $g$ such that $\phi(U_g)$ is Euclidean-open in $V$. Since 
$G = \bigcup_{g \in G} U_g$ and $\phi(G) = \bigcup_{g \in G} \phi(U_g)$, we can conclude that $\rho(G)v$ is Euclidean-open in $V$.
\hfill $\Box$
\\[-1mm]  

\remark
Let $(G, \rho, V)$ be an algebraic triplet over $\C$. 
Then $(G, \rho, V)$ can be regarded as the holomorphic triplet (cf. [5, Ch.1, Sect. 6]), and $(G, \rho, V)$ of each  
category induces the same infinitesimal form. Hence $(G, \rho, V)$ is a PV in the holomorphic category if and only if  
$(G, \rho, V)$ is a PV in the algebraic category. Thus the notion of PV in the holomorphic category is an extension of that in the algebraic category over $\C$. \\[-2mm]  

From now on we consider PVs in the infinitesimal category, and unless otherwise stated  
we always assume that Lie algebars and their representations are defined over the complex number field $\C$. 
Let $(\g, f, V)$ be a triplet composed of a Lie algebra $\g$ and its representation $f$ on $V$.       
In this paper we say that $(\g, f, V)$ is an (infinitesimal) PV if there exists $v \in V$ such that $f(\g)v=V$ 
in view of Proposition \ref{prop1}.   
We also call such an element $v$ a generic point. 
Note that all generic elements form a Euclidean-open set of $V$.   
We say that $(\g, f, V)$ is algebraic if it is a differential of some algebraic triplet $(G, F, V)$ over $\C$.
Then note that an algebraic triplet $(\g, f, V)$ is an (infinitesimal) PV in our sense 
if and only if $(G, F, V)$ is a PV in the algebraic category. 
A triplet $(\g, f, V)$ is said to be irreducible if $f$ is irreducible, and said to be faithful if $f$ is faithful.  
In the following we denote the general linear Lie algebra by $\gl(n)$  instead of $\gl(\C^n)$. 
When $\h$ is a subalgebra of $\gl(n)$,  
we denote the identity representation $\h \to \gl(n)$ of $\h$ by $\Lambda_1$, and denote its dual representation by $\Lambda_1^*$.  

In Theorem \ref{IFPS and hom} we showed that a complex IFPS on $L$ corresponds to an equivalence class $[(f, v)]$, 
where $f: \l \to \gl(n+1)$ is a Lie algebra 
representation and satisfies $f(\l) v \oplus \langle v \rangle = \C^{n+1}$.  Let $\Lambda_1$ be the identity representation of $\gl(1)$, and   
extend the representation $f$ to the tensor product $f \otimes \Lambda_1$.  
Since the tensor product $\C^{n+1} \otimes \C$ is linearly isomorphic to $\C^{n+1}$, we obtain the extended representation   
$f \otimes \Lambda_1:$ $\l \oplus \gl(1) \to \gl(n+1)$.  
Then we have $f \otimes \Lambda_1 (\l \oplus \gl(1)) v = \C^{n+1}$ and therefore $(\l \oplus \gl(1), f \otimes \Lambda_1, \C^{n+1})$ 
is a PV. 
Note that for any non-trivial representation $\alpha: \gl(1) \to \gl(1)$, the tensor product $f \otimes \alpha$ also gives a PV.    
Here we introduce the following notion:  
\begin{df}\label{PV of type IFPS}
Let $\alpha$ be a non-trivial representation $\gl(1) \to \gl(1)$, and let $(\l \oplus \gl(1), f \otimes \alpha, V)$ be a PV. 
Then 
Then we say that $(\l \oplus \gl(1), f \otimes \alpha, V)$ is a PV of type IFPS if it satisfies $\dim \l + 1 = \dim V$.  
\end{df}
\remark
Let $\alpha$ and $\beta$ be non-trivial representations $\gl(1) \to \gl(1)$.  
Then note that the triplet 
$(\l \oplus \gl(1), f \otimes \alpha, V)$ is a PV of type IFPS if and only if 
$(\l \oplus \gl(1), f \otimes \beta, V)$ is a PV of type IFPS. 
In the following we frequently use the identity representation $\Lambda_1$ as a non-trivial representation of $\gl(1)$.    
By the definition if a PV $(\l \oplus \gl(1), f \otimes \alpha, V)$ is of type IFPS, then $f \otimes \alpha$ is faithful. \\[-3mm]

By the above consideration, if $(f, v)$ satisfies $f(\l) v \oplus \langle v \rangle = \C^{n+1}$, then $(\l \oplus \gl(1), f \otimes \Lambda_1, V)$ is 
a PV of type IFPS and $v$ is a generic point. Conversely if $(\l \oplus \gl(1), f \otimes \Lambda_1, V)$ is a PV of type IFPS and $v$ is its generic point, then 
$(f, v)$ satisfies $f(\l) v \oplus \langle v \rangle = \C^{n+1}$.  
Therefore from Theorem \ref{IFPS and hom} we have the following corollary: 
Let $L$ be a complex Lie group, and let $\l$ be its Lie algebra. 
\begin{cor}\label{IFPS and PV} 
There is a one-to-one correspondence between
the set $\{$ complex IFPS on $L$ $\}$ and the set $\{$ $(f \otimes \Lambda_1, v)$ $\mid$ $(\l \oplus \gl(1), f \otimes \Lambda_1, V)$ is a PV of type IFPS, 
$v$ is a generic point of $V$ $\}/\sim$. 
\end{cor}
\remark 
A representation $f \otimes \Lambda_1$ is irreducible if and only if $f$ is irreducible.  Hence an irreducible complex IFPS corresponds to 
an equivalence class $[(f \otimes \Lambda_1, v)]$ such that $f \otimes \Lambda_1$ is irreducible.  

We can consider PVs also in the real $C^\infty$ category, and we have the same assertions from Proposition \ref{prop1} to Corollary 
\ref{IFPS and PV}.  \\[-2mm]

\example
Let us consider a triplet (a): $(\gl(2), 3\Lambda_1, V)$, where $3\Lambda_1$ is the $3$-symmetric product of $\Lambda_1$, and therefore 
we have $\dim V = 4$ (see \cite[p.47]{sato-kimura} about this example). 
Then $3\Lambda_1$ is expressed by matrices as follows: 
{\small
\[
3\Lambda_1 (
\begin{pmatrix}
\alpha & \beta \\
\gamma & \delta
\end{pmatrix}
) 
= 
\begin{pmatrix}
3\alpha & \beta & 0 & 0 \\
3\gamma & 2\alpha + \delta & 2\beta & 0 \\
0 & 2\gamma & \alpha + 2\delta & 3\beta \\
0 & 0 & \gamma & 3\delta
\end{pmatrix}. 
\]
}
By taking $(1,0,0,1)$ as a generic point, it is easily seen that the triplet $(a)$ is a PV. 
Moreover both $\gl(2)$ and $3\Lambda_1$ are decomposed into $\sll(2) \oplus \gl(1)$ and $3\Lambda_1 \otimes \Lambda_1$, and      
$\dim \gl(2) = 4 = \dim V$.  Therefore the triplet (a): ($\sll(2) \oplus \gl(1)$, $3\Lambda_1 \otimes \Lambda_1$, $V(4)$) is of type IFPS. This fact implies that $SL(2, \C)$ admits an (irreducible) complex IFPS. 

\section{Some notions on Prehomogeneous vector space}\label{sec:notion of PV}
In Theorem \ref{main thm} we classify complex Lie groups admitting an irreducible complex IFPS. To prove Theorem \ref{main thm}  
it is sufficient to classify irreducible PVs of type IFPS from Corollary \ref{IFPS and PV}.  
In this section we prepare some notions on PV useful for classifications of PVs. 
Especially we introduce two notions of castling transform and isomorphism, which were already considered in the algebraic   
category in \cite{sato-kimura} by using the terminology of algebraic groups. We prove some properties on these concepts 
needed later for our classification.  
\\[-1mm]

Sato (\cite[p. 37]{sato-kimura}) proved the following proposition in the algebraic category. We denote by $V(n)$ an 
$n$-dimensional vector space, and denote by $\Lambda_1$ the identity representation of $\gl(n)$. 
\begin{prop}\label{cast} 
Let $f$ be a representation of a complex Lie algebra $\h$ on $V(m)$.  
For any $n$ satisfying $m>n\geq 1$, the following conditions are equivalent.
\begin{enumerate}
\item $(\h \oplus \gl(n), f \otimes \Lambda_1, V(m) \otimes V(n))$ is a PV.
\item $(\h \oplus \gl(m-n), f^* \otimes \Lambda_1, V(m)^* \otimes V(m-n))$ is a PV, where $f^*$ is the dual representation of 
$f$.
\end{enumerate}
Generic isotropy subalgebras of these PVs have the same dimensions.
\end{prop}
\pf
The idea of the proof is based on that of Proposition 7 in \cite{sato-kimura}. 
First we deduce (2) from (1). Suppose that $(\h \oplus \gl(n), f \otimes \Lambda_1, V(m) \otimes V(n))$ is a PV.
Then there exists $w \in V(m) \otimes V(n)$ such that 
\[f \otimes \Lambda_1(\h \oplus \gl(n)) w = V(m) \otimes V(n).\] 
We identify $V(m) \otimes V(n)$ and $\underbrace{V(m) \oplus \cdots \oplus V(m)}_n$. 
Then we note that  \[f \otimes \Lambda_1 (H, A) \ x = f(H) \ x + x \ {}^t \! A\]
for $x = (v_1, \cdots, v_n )$ $\in$ $V(m) \otimes V(n)$, $H \in \h$, $A \in \gl(n)$. 
We fix a basis $\{ e_1, \ldots, e_n, e_{n+1}, \ldots, e_m \}$ of $V(m)$ such that $w$ is expressed as $(e_1, \cdots , e_n)$.  
We denote its dual basis by $\{ e_1^*, \ldots, e_n^*, e_{n+1}^*, \ldots, e_m^* \}$.
When we set $w^\perp := (e_{n+1}^*, \cdots , e_m^*)$,  
$w^\perp$ gives a vector of $V(m)^* \otimes V(m-n)$. 
From the assumption we have   
\begin{eqnarray*}
&\ & f \otimes \Lambda_1(\h \oplus \gl(n)) w 
= 
f(\h) 
\begin{pmatrix}
I_n \\
\hline 
0
\end{pmatrix}
+
\begin{pmatrix}
I_n \\
\hline 
0
\end{pmatrix}
\gl(n)  
= 
\begin{pmatrix}
\gl(n) \\
\hline 
M(m-n, n)
\end{pmatrix},
\end{eqnarray*}
where $M(m-n, n)$ denotes the set of $(m-n) \times n$ matrices. 
Hence $f (\h) =$ 
$\left(
\begin{array}{c|c}
\ast & \ast \\
\hline 
M(m-n, n) & \ast
\end{array}
\right)$, which implies that 
$f^*(\h) =$ 
$\left(
\begin{array}{c|c}
\ast & M(n, m-n) \\
\hline 
\ast & \ast
\end{array}
\right)$.
Therefore we have 
\begin{eqnarray*}
&& f^* \otimes \Lambda_1(\h \oplus \gl(m-n))w^\perp \\
&=& 
\left(
\begin{array}{c|c}
\ast & M(n, m-n) \\
\hline 
\ast & \ast
\end{array}
\right)  
\begin{pmatrix}
0 \\
\hline 
I_{m-n}
\end{pmatrix}
+
\begin{pmatrix}
0 \\
\hline 
I_{m-n}
\end{pmatrix}
\gl(m-n)\\
&=&
\begin{pmatrix}
M(n, m-n) \\
\hline 
\ast
\end{pmatrix}
+
\begin{pmatrix}
0 \\
\hline 
\gl(m-n)
\end{pmatrix}\\
&=& V(m)^* \otimes V(m-n).
\end{eqnarray*}

Secondly to deduce (1) from (2), it suffices to apply the result we have just obtained to this case. 

Next we show that $(\h \oplus \gl(n))_w$ is isomorphic to $(\h \oplus \gl(m-n))_{w^\perp}$.
Let $(H, A)$ be an element of $\h \oplus \gl(n)$ such that $f(H)$ is of the form $f(H) =$ 
$\left(
\begin{array}{c|c}
B & C \\
\hline 
D & E
\end{array}
\right)$ .
Then 
$
f \otimes \Lambda_1(H, A) \ w = f(H)\ w + w \ {}^tA =
\begin{pmatrix}
B\\
\hline 
D
\end{pmatrix}
+
\begin{pmatrix}
{}^tA \\
\hline 
0
\end{pmatrix}.$
Hence $(H, A) \in (\h \oplus \gl(n))_w$ if and only if $D=0$, $B+{}^tA=0$. 
Thus we have 
\begin{eqnarray*}
&& (\h \oplus \gl(n))_w \\
&=&
\{ (H, A) \in \h \oplus \gl(n) \mid f(H) = 
\left(
\begin{array}{c|c}
-{}^tA & \ast \\
\hline 
0 & \ast
\end{array}
\right)\} \\
&\cong &  
\{ H \in \h \mid f(H) = 
\left(
\begin{array}{c|c}
\ast & \ast \\
\hline 
0 & \ast
\end{array}
\right) 
\}.
\end{eqnarray*}
Likewise
\begin{eqnarray*}
&& (\h \oplus \gl(m-n))_{w^\perp} \\
&=&
\{ (H, A') \in \h \oplus \gl(m-n) \mid f^*(H) = 
\left(
\begin{array}{c|c}
\ast & 0 \\
\hline 
\ast & -{}^tA'
\end{array}
\right)
\} \\
&=&
\{ (H, A') \in \h \oplus \gl(m-n) \mid f(H) = 
\left(
\begin{array}{c|c}
\ast & \ast \\
\hline 
0 & A'
\end{array}
\right)
\} \\
&\cong &
\{ H \in \h \mid f(H) = 
\left(
\begin{array}{c|c}
\ast & \ast \\
\hline 
0 & \ast
\end{array}
\right)
\}.
\end{eqnarray*}
Therefore $(\h \oplus \gl(n))_w \cong (\h \oplus \gl(m-n))_{w^\perp}$.
\hfill $\Box$
\\[-1mm]

We call the transformation $(1) \longleftrightarrow (2)$ a c-transformation in this paper. 
We can easily verify that c-transformations preserve 
irreducibility by using the following fact: 
let $(\g_1, f_1, V_1)$ and $(\g_2, f_2, V_2)$ be irreducible triplets. Then the tensor product 
$f_1 \otimes f_2$ is irreducible if and only if both $f_1$ and $f_2$ are irreducible.   
 
Now let us consider the triplet 
$(\h \oplus \gl(n), f \otimes \Lambda_1, V(m) \otimes V(n))$ which satisfies $\dim \h + n^2 = mn$.  
This triplet is equal to the the triplet 
$(\h \oplus \sll(n) \oplus \gl(1), f \otimes \Lambda_1 \otimes \Lambda_1, V(m) \otimes V(n) \otimes V(1))$.  
In this context if this triplet is of type IFPS, then 
the c-transformation of this triplet gives again a PV of type IFPS.  Indeed  
\[(\h \oplus \gl(n), f \otimes \Lambda_1, V(m) \otimes V(n))\]
is c-transformed into the triplet 
\[(\h \oplus \gl(m-n), f^* \otimes \Lambda_1, V(m)^* \otimes V(m-n)),\] 
which is equal to 
\[(\h \oplus \sll(m-n) \oplus \gl(1), f^* \otimes \Lambda_1 \otimes \Lambda_1, V(m)^* \otimes V(m-n) \otimes V(1)).\] 
Since we have $\dim \h + (m-n)^2 = m(m-n)$,  
it follows that  the last triplet is of type IFPS. \\[-3mm] 

\example
Let us illustrate the c-transformations by using the example (a) in \S \ref{sec:PV and IFPS}.   
The triplet (a): $(\gl(2), 3\Lambda_1, V(4))$ is naturally identified with  
$(\sll(2) \oplus \gl(1), 3\Lambda_1 \otimes \Lambda_1, V(4) \otimes V(1))$.  Hence by a c-transformation of the triplet (a) we obtain the triplet 
$(\sll(2) \oplus \gl(3), 3\Lambda_1^* \otimes \Lambda_1, V(4)^* \otimes V(3))$. This triplet is equal to 
$(\sll(2) \oplus \sll(3) \oplus \gl(1), 3\Lambda_1^* \otimes \Lambda_1 \otimes \Lambda_1, V(4)^* \otimes V(3) \otimes V(1))$, which is of 
type IFPS. 

\begin{df}[{\cite[p.245]{kimura}}]\label{def of isom}
Two triplets $(\g, \rho, V)$ and $(\g', \rho', V')$ {\rm(}not necessarily PVs{\rm)} are said to be isomorphic if there exists a Lie algebra isomorphism 
$\sigma$: $\rho(\g) \to \rho'(\g')$ and a linear isomorphism $\tau: V \to V'$ satisfying 
$\sigma \circ \rho (X)$ $=$ $\tau \rho(X) \tau^{-1}$ for any $X \in \g$.  
Then we write  $(\g, \rho, V)$ $\cong$ $(\g', \rho', V')$. 
\end{df} 
Note that when $(\g, \rho, V)$ and $(\g', \rho', V')$ are isomorphic, they are said to be strongly equivalent in \cite[p.36]{sato-kimura}. 
Moreover in \cite[p.36]{sato-kimura}, the isomorphism of triplets is defined by using the terminology of algebraic groups as follows: 
let ($G$, $\rho$, $V$) and ($G'$, $\rho'$, $V'$) be algebraic triplets. Then we denote $(G, \rho, V)$ $\cong$ $(G', \rho', V')$ if there exists 
a rational isomorphism $\sigma$: $\rho(G) \to \rho'(G')$ and a linear isomorphism $\tau: V \to V'$ satisfying 
$\sigma \circ \rho (g)$ $=$ $\tau \rho(g) \tau^{-1}$ for any $g \in G$. 
Clearly by differentiating $\sigma$, we obtain the isomorphism $d\sigma$: $(Lie(G), d\rho, V)$ $\to$ $(Lie(G'), d\rho', V')$, and we have 
$(Lie(G), d\rho, V)$ $\cong$ $(Lie(G'), d\rho', V')$ in the above sense.  

For example, let $(\g, \rho, V)$ be an arbitrary triplet. We denote  the natural inclusion $\rho(\g) \hookrightarrow \gl(V)$ by $\imath$. 
Then obviously we have $(\g, \rho, V)$ $\cong$ $(\rho(\g), \imath, V)$. 
 
Note that the isomorphism $\cong$ is an equivalence relation. We call the equivalence class of $(\g, \rho, V)$
an isomorphism class  of $(\g, \rho, V)$, and denote it by $[(\g, \rho, V)]$.  
We here introduce a new concept. Let $\h$ and $\h'$ be subalgebras of $\g$ and $\g'$ respectively. We say that $\rho(\h)$ corresponds to $\rho'(\h')$ (via $\sigma$) if $\sigma(\rho(\h)) = \rho'(\h')$.  
\\[-2mm] 

\remark 
Consider two isomorphic triplets $(\g, \rho, V)$ $\cong$ $(\g', \rho', V')$.  Then we have the following: \\
(1) If $(\g, \rho, V)$ is a PV, then $(\g', \rho', V')$ is also a PV. \\ 
(2) If $(\g, \rho, V)$ is irreducible, then $(\g', \rho', V')$ is also irreducible. \\
However note that $\g$ and $\g'$ are not necessarily isomorphic as Lie algebras.  

We use the following proposition in \S 7. 
\begin{prop}\label{isom of algebraic}
Let $(G, \rho ,V)$ and $(G', \rho', V')$ be algebraic triplets, and let $(\g, d\rho, V)$ and $(\g', d\rho', V')$ be the induced triplets 
respectively. If  $(\g, d\rho, V)$ $\cong$ $(\g', d\rho', V')$,  
then  we have $(G, \rho ,V)$ $\cong$ $(G', \rho', V')$. 
\end{prop}
\pf 
Since $\rho$ is a rational representation and $G$ is a connected algebraic group, its image $\rho(G)$ is also a connected linear algebraic 
subgroup of $GL(V)$ (cf. \cite[p.102]{onishchik}). Thus $\rho'(G')$ is also one of $GL(V')$. 
From the assumption, there exists a Lie algebra isomorphism $\sigma$: $d\rho(\g) \to d\rho'(\g')$ and a linear isomorphism 
$\tau$: $V \to V'$ such that 
$\sigma (d\rho(X)) = \tau d\rho(X) \tau^{-1}$ for $X \in \g$. 

Now we define the map $\tilde{\sigma}$: $\rho(G) \to GL(V')$ by $\tilde{\sigma}(\rho(g)) = \tau \rho(g) \tau^{-1}$ for $g \in G$. 
Then $\tilde{\sigma}$ is a rational isomorphism, and its differential $d\tilde{\sigma}$ is equal to $\sigma$. 
Thus the image $\tilde{\sigma}(\rho(G))$ is a connected linear algebraic subgroup of $GL(V')$.   
Moreover its Lie algebra is $\sigma(d\rho(\g)) = d\rho'(\g')$, it follows that $\tilde{\sigma}(\rho(G))$ $=$ $\rho'(G')$. 
Hence $(G, \rho ,V)$ $\cong$ $(G', \rho', V')$.   
\hfill $\Box$\\[-2mm]

Here we prove some elementary properties on isomorphism of triplets we shall use later.  
\begin{prop}\label{prop of isom}
\ 
\begin{enumerate}
\item Suppose that $(\g, f, V) \cong (\g', f', V')$ and $(\h, g, W) \cong (\h', g', W')$.  
Then $(\g \oplus \h, f \otimes g, V \otimes W)$ $\cong$ $(\g' \oplus \h', f' \otimes g', V' \otimes W')$. 
\item Let $(\g \oplus \h, f \otimes g, V \otimes W)$ and $(\g' \oplus \h', f' \otimes g', V' \otimes W')$ be irreducible triplets. 
Suppose that these triplets are isomorphic and $f \otimes g(\g)$ corresponds to $f' \otimes g'(\g')$. 
Then $(\g, f, V)$ $\cong$ $(\g', f', V')$.  
\item Consider faithful triplets $(\g, \rho, V)$ and $(\g', \rho', V')$.    
If these triplets are isomorphic, then $\g$ is isomorphic to $\g'$ as a Lie algebra.  
\item Let $\alpha$ and $\beta$ be non-trivial representations $\gl(1) \to \gl(1)$. Then   
$(\gl(1), \alpha, V(1))$ is isomorphic to $(\gl(1), \beta, V(1))$. 
\item 
Let $\h$ be a semisimple Lie algebra. Then any triplet $(\h, \rho, V)$ is isomorphic to its dual $(\h, \rho^*, V^*)$. 
\end{enumerate} 
\end{prop}
\pf 
(1) 
From the assumption there exist isomorphisms $\phi: f(\g) \to f'(\g')$, $\psi: g(\h) \to g'(\h')$, and linear isomorphisms $\tau: V \to V'$, 
$\upsilon: W \to W'$ such that $\phi\circ f(X) = \tau f(X) \tau^{-1}$ and $\psi \circ g(Y) = \upsilon g(Y) \upsilon^{-1}$ 
for $X \in \g$ and $Y \in \h$.   
We define a map 
\[\phi \otimes \psi: f \otimes g(\g \oplus \h) \to f' \otimes g'(\g' \oplus \h')\]  
by \ $f(X) \otimes I_W + I_V \otimes g(Y)$ $\mapsto$ $\phi \circ f(X) \otimes I_{W'} + I_{V'} \otimes \psi \circ g(Y)$.  
Then this mapping is well defined, and we can easily check that this is a Lie algebra isomorphism. 
Moreover for $X, Y \in \g \oplus \h$, 
\begin{eqnarray*}
\phi \otimes \psi  \circ f \otimes g(X, Y) &=& \phi \circ f(X) \otimes I_{W'} + I_{V'} \otimes \psi \circ g(Y) \\ 
&=& \tau \otimes \upsilon (f \otimes g(X, Y)) (\tau \otimes \upsilon)^{-1}.  
\end{eqnarray*}
Hence we obtain our claim. \\[-2mm]

(2) 
From the assumption, we have $(\g \oplus \h, f \otimes g, V \otimes W)$ $\cong$ $(\g' \oplus \h', f' \otimes g', V' \otimes W')$. Hence 
there exist an isomorphism $\sigma: f \otimes g(\g \oplus \h) \to f' \otimes g'(\g' \oplus \h')$ 
and a linear isomorphism $\tau: V \otimes W \to V' \otimes W'$ such that 
\[\sigma(f(X) \otimes I_W + I_V \otimes g(Y)) \tau(x) = \tau ((f(X) \otimes I_W + I_V \otimes g(Y))(x))\] for $(X, Y) \in \g \oplus \h$ and 
$x \in V \otimes W$. 
Since $f \otimes g(\g)$ corresponds to $f' \otimes g'(\g')$, 
$\sigma$ is restricted to the isomorphism $f(\g) \otimes I_W$ $\to$ $f'(\g') \otimes I_{W'}$.   
Let $\imath$ be the natural isomorphism $f(\g) \to f(\g) \otimes I_W$ defined by $f(X) \mapsto f(X) \otimes I_W$ for $X \in \g$, 
and let $\imath'$ be the natural isomorphism $f'(\g') \to f'(\g') \otimes I_{W'}$ defined by the same way. 
We denote by the same symbol $\sigma$ the composite of isomorphisms 
\[f(\g) \to f(\g) \otimes I_W \to f'(\g') \otimes I_{W'} \to f'(\g'),\] which is equal to 
$\imath'^{-1} \circ \sigma \circ \imath$.   

Let $\{e_1, \cdots, e_n\}$ be a basis of $W$. 
Since $(\g, f, V)$ is irreducible,  $(\g, f \otimes I_W, V \otimes e_1)$ is an irreducible $\g$-submodule of $V \otimes W$.  
Then we have the natural $\g$-isomorphism $\phi: V \to V \otimes e_1$ defined by $\phi(v) = v \otimes e_1$.     
Since $f(\g) \otimes I_W$ is isomorphic to $f'(\g') \otimes I_{W'}$ via $\sigma$ and we have 
$\sigma(f(X) \otimes I_W ) \tau(v \otimes e_1)$ $=$ $\tau (f(X) \otimes I_W (v \otimes e_1))$ for $X \in \g$ and $v \in V$,  
it follows that $(\g', f' \otimes I_{W'}, \tau(V \otimes e_1))$ 
is also an irreducible $\g'$-submodule of $V' \otimes W'$. 
Let $\{e'_1, \cdots, e'_{n'}\}$ be a basis of $W'$. 
Then $V' \otimes W'$ is decomposed into the direct sum of equivalent irreducible $\g'$-submodules 
$V' \otimes e'_1 \oplus \cdots \oplus V' \otimes e'_{n'}$.
We note that $f'(\g')$ is reductive with at most one-dimensional center $\{\lambda I_{V'} \mid \lambda \in \C\}$  
because $f'$ is irreducible (cf. \cite[p.2]{sato-kimura}).  
Thus concerning the triplet $(\g', f' \otimes I_{W'}, V' \otimes W')$, its image $f' \otimes I_{W'}(\g') = f'(\g') \otimes I_{W'}$ is reductive 
whose center is contained in $\{\lambda I_{V' \otimes W'} \mid \lambda \in \C\}$. 
It follows that 
($\g'$, $f' \otimes I_{W'}$, $V' \otimes W'$) is completely reducible.  
Since $(\g', f' \otimes I_{W'}, \tau(V \otimes e_1))$ is an irreducible $\g'$-submodule of $V' \otimes W'$, 
there exists a $\g'$-isomorphism $\tau(V \otimes e_1) \to V' \otimes e'_1$ since 
$V' \otimes W'$ is the direct sum of equivalent irreducible $\g'$-submodules $V' \otimes e'_1 \oplus \cdots \oplus V' \otimes e'_{n'}$  
as stated before.  
   
This gives the sequence of $\g'$-isomorphisms  
\[\tau(V \otimes e_1) \to V' \otimes e'_1 \to V', \] 
and we denote the composite of these by $\psi$. By using $\psi$ let us consider the composite  
$\psi \circ \tau \circ \phi$, which gives a linear isomorphism $V \to V'$.  
Then for $X \in \g$ and $v \in V$, we  have 
\begin{eqnarray*}
\sigma(f(X)) \psi \circ \tau \circ \phi(v) &=& \psi(\sigma(f(X) \otimes I_W) \tau \circ \phi(v)) \\  
&=& \psi \circ \tau(f(X) \otimes I_W (\phi(v))) \\
&=& \psi \circ \tau \circ \phi(f(X)v).
\end{eqnarray*}
Hence $(\g, f, V) \cong (\g', f', V')$.  \\[-2mm]

(3) 
From the assumption there exists a Lie algebra isomorphism $\sigma: \rho(\g) \to \rho'(\g')$. We set 
$\phi := \rho'^{-1} \circ \sigma \circ \rho$, then $
\phi$ gives a Lie algebra isomorphism from $\g$ to $\g'$. \\[-2mm]

(4) Since $\alpha(\gl(1))$ $=$ $\beta(\gl(1))$ $=$ $\gl(1)$, immediately we obtain our claim. \\[-2mm]

(5) 
The triplet $(\h, \rho, V)$ is isomorphic to the algebraic triplet $(\rho(\h), \imath, V)$. We prove this fact in Proposition \ref{algebraic},  
which we explain later. 
Thus there exits an algebraic triplet 
$(\tilde{H}, \tilde{\imath}, V)$ such that its infinitesimal form is $(\rho(\h), \imath, V)$. Then it has been proved that 
$(\tilde{H}, \tilde{\imath}, V)$ is isomorphic to its dual $(\tilde{H}, \tilde{\imath}^*, V^*)$ (cf. \cite[p.245]{kimura}). 
Since the infinitesimal form of $(\tilde{H}, \tilde{\imath}^*, V^*)$ is equal to $(\rho(\h), \imath^*, V^*)$, 
by differentiation we have $(\rho(\h), \imath, V)$ $\cong$ $(\rho(\h), \imath^*, V^*)$. 
Moreover this triplet $(\rho(\h), \imath^*, V^*)$ is isomorphic to 
$(\h, \rho^*, V^*)$. Hence we have $(\h, \rho, V)$ $\cong$ $(\h, \rho^*, V^*)$. 
\hfill $\Box$\\[-2mm]

By combining two notions isomorphism and c-transformation,  a castling transform is defined in \cite{sato-kimura}. This 
is an important tool for classification of irreducible PVs.  
\begin{df}[{\cite[p.39]{sato-kimura}}]\label{def of cast}
Two triplets $(\g, \rho, V)$ and $(\g', \rho', V')$ {\rm(}not necessarily PVs{\rm)} are said to be castling transforms of each other when there exist a triplet 
$(\tilde{\g}, \tilde{\rho}, V(m))$ and a positive number $n$ with $m > n \geq 1$ such that 
\[(\g, \rho, V) \cong (\tilde{\g} \oplus \sll(n), \tilde{\rho} \otimes \Lambda_1, V(m) \otimes V(n))\] 
and 
\[(\g', \rho', V') \cong (\tilde{\g} \oplus \sll(m-n), \tilde{\rho}^* \otimes \Lambda_1, V(m)^* \otimes V(m-n)).\]  
A triplet $(\g, \rho, V)$ is said to be reduced if there is no castling transform $(\g', \rho', V')$ of $(\g, \rho, V)$ with  
$\dim V' < \dim V$. 
\end{df}
Note that if two triplets $(\g, \rho, V)$ and $(\g', \rho', V')$ are isomorphic, then 
$(\g, \rho, V)$ is reduced if and only if  $(\g', \rho', V')$ is reduced. 

Similarly to the case of isomorphism $\cong$, originally castling transform is defined by using the terminology of algebraic groups 
(see \cite[p.39]{sato-kimura}).  
Concerning Definition \ref{def of cast}, let us consider the transformations:
\begin{eqnarray*}
&&(1) \ (\tilde{\g} \oplus \sll(n), \tilde{\rho} \otimes \Lambda_1, V(m) \otimes V(n)) \\
&&(2) \ (\tilde{\g} \oplus \sll(m-n), \tilde{\rho}^* \otimes \Lambda_1, V(m)^* \otimes V(m-n)),
\end{eqnarray*}
and
\begin{eqnarray*} 
&&(1)^* \ (\tilde{\g} \oplus \sll(n), \tilde{\rho} \otimes \Lambda_1^*, V(m) \otimes V(n)^*) \\  
&&(2)^* \ (\tilde{\g} \oplus \sll(m-n), \tilde{\rho}^* \otimes \Lambda_1^*, V(m)^* \otimes V(m-n)^*). 
\end{eqnarray*} 
We call each transformation $(1) \longleftrightarrow (2)$, $(1)^* \longleftrightarrow (2)^*$ an sc-transformation.  
Now concerning a triplet $(\tilde{\g} \oplus \sll(n), \tilde{\rho} \otimes \Lambda_1, V(m) \otimes V(n))$,  
suppose that there exist a Lie algebra $\l$ and its representation $f: \l \to \gl(V(m))$ such that $\tilde{\g}$ = $\l \oplus \gl(1)$ and 
$\tilde{\rho} = f \otimes \alpha$, where $\alpha$ is a non-trivial representation $\gl(1) \to \gl(1)$.  
We suppose that this triplet is a PV.     
Then from (1) and (4) of Proposition \ref{prop of isom},  
\[(\l \oplus \gl(1) \oplus \sll(n), f \otimes \alpha \otimes \Lambda_1, V(m)\otimes V(1) \otimes V(n))\] 
is isomorphic to  
\[(\l \oplus \gl(1) \oplus \sll(n), f \otimes \Lambda_1 \otimes \Lambda_1, V(m)\otimes V(1) \otimes V(n)).\]  
It follows that   
by Proposition \ref{cast} 
\[(1) \ (\tilde{\g} \oplus \sll(n), \tilde{\rho} \otimes \Lambda_1, V(m) \otimes V(n)) \quad  \mbox{is a PV if and only if} \] 
\[\hspace{-4mm} (2) \ (\tilde{\g} \oplus \sll(m-n), \tilde{\rho}^* \otimes \Lambda_1, V(m)^* \otimes V(m-n)) \ \ \mbox{is a PV}.\] 
Moreover by the same way we can show that $(1)^*$: $\tilde{\rho} \otimes \Lambda_1^*$ is a PV if and only if $(2)^*$: $\tilde{\rho}^* \otimes \Lambda_1^*$
is a PV.  Hence as before  
we see that a PV of type IFPS is transformed into another PV of type IFPS by any sc-transformation.   
\begin{df}[{\cite[p.39]{sato-kimura}}] 
Two triplets $(\g, \rho, V)$ and $(\g', \rho', V')$ are said to be castling equivalent when one is obtained from the other by a finite number of 
castling transforms. We call this equivalence class a castling class. 
\end{df} 
For any triplet $(\g, \rho, V)$, since $\dim V$ is finite, its castling class has at least one reduced triplet. 
Furthermore we can prove that the castling class of an irreducible triplet contains only one reduced triplet up to isomorphism.
For the proof see \cite[p.39]{sato-kimura}. By applying the same argument of \cite[p.39]{sato-kimura} to the infinitesimal category, 
we obtain the assertion. Actually we do not use this result, and we omit its proof here.  

We use the following assertions in \S 7 to determine all triplets castling equivalent to a reduced one. 

\begin{prop}\label{prop iso and cast}
Let $(\g, \rho ,V)$ be a faithful irreducible triplet. Then for any castling transform $(\g', \rho' ,V')$ of $(\g, \rho ,V)$, there exists 
an sc-transformation of $(\g, \rho ,V)$ isomorphic to $(\g', \rho' ,V')$. 
\end{prop}
\pf 
From the assumption there exists a triplet $(\tilde{\g}, \tilde{\rho}, V(m))$ such that
\[(\g, \rho ,V) \cong (\tilde{\g} \oplus \sll(n), \tilde{\rho} \otimes \Lambda_1, V(m) \otimes V(n)) \ \ \mbox{and}\]
\[(\g', \rho' ,V') \cong (\tilde{\g} \oplus \sll(m-n), \tilde{\rho}^* \otimes \Lambda_1, V(m)^* \otimes V(m-n)). \] 
Then by the definition of isomorphism of triplets, there exist an isomorphism $\sigma:$ $\tilde{\rho} \otimes \Lambda_1 (\tilde{\g} \oplus \sll(n))$ $\to$ $\rho(\g)$ 
and $\tau: V(m) \otimes V(n) \to V$ such that 
\[\sigma(\tilde{\rho}(X) \otimes I_n + I_m \otimes Y) \tau (v) = \tau ((\tilde{\rho}(X) \otimes I_n + I_m \otimes Y) (v)) \] 
for $(X, Y) \in \tilde{\g} \oplus \sll(n)$ and $v \in V(m) \otimes V(n)$.  
Now we claim that 
there exists a (reductive) ideal $\h$ of $\g$ such that 
$\g = \h \oplus \sll(n)$ and $\sigma(\tilde{\rho} \otimes \Lambda_1(\tilde{\g}))$ $=$ $\rho(\h)$ and 
$\sigma(\tilde{\rho} \otimes \Lambda_1(\sll(n)))$ $=$ $\rho(\sll(n))$.
First since $\rho$ is faithful, we have the sequence of isomorphisms of Lie algebras  
\[\g \to \rho(\g) \to \tilde{\rho} \otimes \Lambda_1 (\tilde{\g} \oplus \sll(n))\] 
given by $\sigma^{-1} \circ \rho$.  
Next we observe that $\tilde{\rho} \otimes \Lambda_1 (\tilde{\g} \oplus \sll(n))$ is equal to the direct sum of Lie algebras 
$\tilde{\rho}(\tilde{\g}) \otimes I_n \oplus I_m \otimes \sll(n)$.  
Therefore by putting $\h$ $=$ $(\sigma^{-1} \circ \rho)^{-1} (\tilde{\rho}(\tilde{\g}) \otimes I_n)$, we obtain a decomposition  
$\g = \h \oplus \sll(n)$. Moreover we can easily see that $\sigma(\tilde{\rho} \otimes \Lambda_1(\tilde{\g}))$ $=$ $\rho(\h)$ and  
$\sigma(\tilde{\rho} \otimes \Lambda_1(\sll(n)))$ $=$ $\rho(\sll(n))$. 
Then there exist irreducible representations $f: \h \to \gl(V(k))$ and $g: \sll(n) \to \gl(V(l))$ such that 
$\rho = f \otimes g$ (cf. \cite[p.236]{kimura}). Hence $(\g, \rho ,V)$ $=$ $(\h \oplus \sll(n), f \otimes g, V(k) \otimes V(l))$. 
Note that $\sigma$ is restricted to the isomorphisms  $\tilde{\rho} \otimes \Lambda_1(\tilde{\g})$ $\to$ $f \otimes g(\h)$ and 
$\tilde{\rho} \otimes \Lambda_1(\sll(n))$ $\to$ $f \otimes g(\sll(n))$.  
Then by (2) in Proposition \ref{prop of isom}, we have  $(\tilde{\g}, \tilde{\rho}, V(m))$ $\cong$ $(\h, f, V(k))$ and 
 $(\sll(n), \Lambda_1, V(n))$ $\cong$ $(\sll(n), g, V(l))$, 
thus especially  $m = k$ and $n = l$.   
Since $g$ is an irreducible representation of $\sll(n)$ with degree $n$, $g$ is equivalent to the identity representation $\Lambda_1$ or its dual $\Lambda_1^*$.  
Therefore $(\g, \rho ,V)$ = $(\h \oplus \sll(n), f \otimes \Lambda_1^{(*)}, V(m) \otimes V(n)^{(*)})$. 
Thus $(\g, \rho ,V)$ can be sc-transformed into the triplet: 
\[(a) \ (\h \oplus \sll(m-n), f^* \otimes \Lambda_1^{(*)}, V(m)^* \otimes V(m-n)^{(*)}).\]  
Since $(\h, f, V(m))$ $\cong$ $(\tilde{\g}, \tilde{\rho}, V(m))$, we have $(\h, f^*, V(m)^*)$ $\cong$ $(\tilde{\g}, \tilde{\rho}^*, V(m)^*)$.  
On the other hand  $(\sll(m-n), \Lambda_1^*, V(m-n)^*)$ $\cong$ 
$(\sll(m-n), \Lambda_1, V(m-n))$ from (5) of Proposition \ref{prop of isom}. Hence  
the triplet $(a)$ is isomorphic to $(\g', \rho' ,V')$ by (1) of Proposition \ref{prop of isom}. 
Therefore the proposition follows.  
\hfill $\Box$ \\

Since sc-transformations preserve faithfulness of triplets, we have the following: 
\begin{cor}\label{cast of reduced triplet}
Let $(\g, \rho ,V)$ be a faithful irreducible triplet. Suppose that a triplet $(\g', \rho' ,V')$ is castling equivalent to 
$(\g, \rho, V)$. Then there exists a triplet obtained by a finite number of sc-transformations from $(\g, \rho, V)$ which is isomorphic to $(\g', \rho' ,V')$. 
\end{cor}

\section{Classification of irreducible PVs of type IFPS} \label{sec:classification}
In this section we shall prove Theorem \ref{main thm}. For this purpose 
it is sufficient to 
classify isomorphism classes of irreducible infinitesimal PVs of type IFPS $(\g, f, V)$ by Corollary \ref{IFPS and PV} and (3) in Proposition \ref{prop of isom}. 
Here Sato and Kimura \cite{sato-kimura} have classified isomorphism classes of reduced irreducible PVs $(G, F, V)$ in the algebraic category. 
In the following we show that  
by differentiating isomorphism classes of reduced irreducible algebraic PVs $(G, F, V)$, we can directly obtain  
a classification of 
isomorphism classes of reduced irreducible infinitesimal PVs (in the holomorphic category).    
First we show the following proposition: 
\begin{prop}\label{algebraic}
\ 
\begin{enumerate}
\item Let $\h$ be a semisimple Lie algebra. 
Then for any triplet $(\h, \rho, V)$,  the triplet $(\rho(\h), \imath, V)$ is algebraic. 
\item Let $\g$ be a Lie algebra, and let $f: \g \to \gl(V)$ be a irreducible representation. 
Then the triplet $(f(\g), \imath, V)$ is algebraic. 
\end{enumerate}
\end{prop}
\pf
(1) First we show that $\rho(\h)$ is algebraic, i.e., there exists a connected semisimple linear algebraic 
group $\tilde{H} \subset GL(V)$ such that $Lie (\tilde{H})$ $=$ $\rho(\h)$. 
To see this, let $H$ be a simply connected complex Lie group with Lie algebra $\h$, and let $\tilde{\rho}$ be the representation $H \to GL(V)$ such that 
$d\tilde{\rho} = \rho$. 
Then $H$ is algebraic and $\tilde{\rho}$ is a rational representation (cf. \cite[p.30-31]{onishchik-vinberg}). 
Therefore $\tilde{\rho}(H)$ is a connected linear algebraic group (cf. \cite[p.102]{onishchik}), and its Lie algebra is $\rho(\h)$. 
Then the inclusion $\tilde{\imath}: \tilde{\rho}(H) \hookrightarrow GL(V)$ is a rational representation of $\tilde{\rho}(H)$, and 
$d\tilde{\imath}$ is equal to the inclusion $\imath$: $\rho(\h) \hookrightarrow \gl(V)$.  
Hence $(\tilde{\rho}(H), \tilde{\imath}, V)$ is an algebraic triplet whose infinitesimal form is 
$(\rho(\h), \imath, V)$. Since this triplet $(\rho(\h), \imath, V)$ is isomorphic to $(\h, \rho, V)$, we obtain our claim. \\[-3mm]

(2) From the assumption $f(\g)$ is a reductive Lie algebra with at most one-dimensional center.  
First we consider the case that $f(\g)$ is semisimple. 
Then from the claim (1), the triplet $(f(\g), \imath, V)$ itself is algebraic. 

Next we consider the case that $f(\g)$ has a one dimensional center $\langle I_V \rangle$.  
Now we denote the semisimple part of $f(\g)$ by $\h$, hence we have $f(\g)$ $=$ $\h \oplus \langle I_V \rangle$. Then from the claim (1) 
there exists a connected linear algebraic group $H$ whose Lie algebra is $\h$. 
Concerning the connected algebraic group 
$H \times GL(1)$, we define the map $\Lambda_1 \otimes \Lambda_1$: $H \times GL(1) \to GL(V)$ by   
$\Lambda_1 \otimes \Lambda_1 (a, \lambda)$ $=$ $\lambda a$. Then this is an algebraic group homomorphism,  
therefore its image $\Lambda_1 \otimes \Lambda_1(H \times GL(1))$ is a connected linear algebraic group of $GL(V)$ (cf. \cite[p.102]{onishchik}).  By differentiating 
$\Lambda_1 \otimes \Lambda_1$, we obtain a Lie algebra representation 
$\Lambda_1 \otimes \Lambda_1$: $\h \oplus \gl(1) \to \gl(V)$. 
Hence the image $\Lambda_1 \otimes \Lambda_1(H \times GL(1))$ has the Lie algebra 
$\Lambda_1 \otimes \Lambda_1(\h \oplus \gl(1))$, and this is equal to $f(\g)$. 
Consequently we obtain the algebraic triplet ($\Lambda_1 \otimes \Lambda_1(H \times GL(1))$, $\imath$, $V$) whose infinitesimal form is equal to 
$(f(\g), \imath, V)$.  \hfill $\Box$ \\

Now let us examine that the set of isomorphism classes of reduced irreducible 
infinitesimal PVs (in the holomorphic category) is given from  
the set of isomorphism classes of reduced irreducible algebraic PVs.      
To see this we define the map  
\begin{eqnarray*} 
&&\Delta: \{\mbox{reduced irreducible algebraic PV $(G, F, V)$}\}/_{\cong} \  \to \\
&& \quad \{\mbox{reduced irreducible infinitesimal PV $(\g, f, V)$}\}/_{\cong}
\end{eqnarray*}
by differentiation, i.e., $\Delta([(G, F, V)]) = [(Lie(G), dF, V)]$. 
Then we can verify that this map is well defined. Here we show that a reduced triplet is mapped into a reduced triplet:   
let $(G, F, V)$ be a reduced irreducible triplet, and we put  $(\g, f, V)$ $:=$ $(Lie(G), dF, V)$. 
We assume that $(\g, f, V)$ is not reduced, thus there exists a triplet 
$(\tilde{\g}, \tilde{\rho}, V(m))$ and a positive number $n$ with $m > n \geq 1$ such that 
\[(\g, f, V) \cong  (\tilde{\g} \oplus \sll(n), \tilde{\rho} \otimes \Lambda_1, V(m) \otimes V(n)) \] 
and $mn > m(m-n)$. 
The last triplet is isomorphic to the triplet
\[(a) \ \ (\tilde{\rho}(\tilde{\g}) \oplus \sll(n), \imath \otimes \Lambda_1, V(m) \otimes V(n)).\]
Since $\tilde{\rho}$ is irreducible, by (2) in Proposition \ref{algebraic}  
($\tilde{\rho}(\tilde{\g})$, $\imath$, $V(m)$) is algebraic, i.e., there exists an algebraic triplet 
($\tilde{G}$, $\tilde{\imath}$, $V$) whose infinitesimal form is  ($\tilde{\rho}(\tilde{\g})$, $\imath$, $V(m)$).  
Now let $\Lambda_1: SL(n) \to GL(n)$ be the identity representation of $SL(n)$, and we   
consider the triplet 
\[(A) \ \ (\tilde{G} \times SL(n), \tilde{\imath} \otimes \Lambda_1, V(m) \otimes V(n)).\] 
Then this is algebraic   
because $\tilde{G} \times SL(n)$ is naturally regarded as a connected linear algebraic subgroup of $GL(m+n)$, and 
$\tilde{\imath} \otimes \Lambda_1$ is a rational representation of $\tilde{G} \times SL(n)$.  
Moreover by differentiating the triplet $(A)$, we obtain the triplet $(a)$.   
Since the triplet $(a)$ is isomorphic to $(\g, f, V)$,   
by Proposition \ref{isom of algebraic} we have $(G, F, V)$ $\cong$ ($\tilde{G} \times SL(n)$, $\tilde{\imath} \otimes \Lambda_1$, $V(m) \otimes V(n)$). 
Hence ($\tilde{G} \times SL(m-n)$, $\tilde{\imath}^* \otimes \Lambda_1$, $V(m)^* \otimes V(m-n)$) is a castling transform of $(G, F, V)$. 
From the hypothesis we have $mn > m(m-n)$, and therefore $(G, F, V)$ is not reduced. This is a contradiction.  
Hence we conclude that $(\g, f, V)$ is reduced. 
\\[-3mm]
 
Now we claim that $\Delta$ is bijective.  
We first show that $\Delta$ is injective. 
Suppose that $\Delta([(G, F, V)])$ $=$ $\Delta([(G', F', V')])$. Then their infinitesimal forms $(Lie(G), dF, V)$ and $(Lie(G'), dF', V')$ 
are isomorphic. Hence by Proposition \ref{isom of algebraic},  
we have $(G, F, V)$ $\cong$ $(G', F', V')$.    
Therefore $\Delta$ is injective. 

Next we show that $\Delta$ is surjective.  
Let $(\g, f, V)$ be a reduced irreducible infinitesimal PV. 
Then the triplet $(\g, f, V)$ is isomorphic to $(f(\g), \imath, V)$, and by (2) in Proposition \ref{algebraic} this triplet is algebraic. 
Hence there exists an algebraic triplet $(\tilde{G}, \tilde{\imath}, V)$ whose infinitesimal form is $(f(\g), \imath, V)$. 
Now suppose that $(\tilde{G}, \tilde{\imath}, V)$ is not reduced, thus there exists a castling transform $(\tilde{G'}, \tilde{\imath}', V')$ of 
$(\tilde{G}, \tilde{\imath}, V)$ such that $\dim V' < \dim V$.  
Then by differentiating these triplets, we obtain a castling transform $(Lie(\tilde{G'}), d\tilde{\imath}', V')$ of $(\g, f, V)$. 
This is a contradiction. Hence $(\tilde{G}, \tilde{\imath}, V)$ is a reduced irreducible algebraic triplet whose infinitesimal form is  
$(f(\g), \imath, V)$. Since  $(\tilde{G}, \tilde{\imath}, V)$ is a PV, $\Delta$ is surjective. 

Hence we conclude that $\Delta$ is bijective.  \\     

Now we begin our classification. Note that any irreducible PV of type IFPS 
is isomorphic to a PV of type IFPS of the form $(\l \oplus \gl(1), \rho \otimes \Lambda_1, V)$ 
because $(\gl(1), \alpha, V(1))$ is isomorphic to $(\gl(1), \Lambda_1, V(1))$ and we have (1) in Proposition \ref{prop of isom}. 
Furthermore its image $\rho \otimes \Lambda_1(\l \oplus \gl(1))$ is a reductive Lie algebra with one dimensional center because 
$\rho \otimes \Lambda_1$ is faithful(cf. the remark after Definition \ref{PV of type IFPS}). Hence we shall investigate only such PVs.   
Here we note that isomorphism classes of reduced irreducible algebraic PVs ($G$, $\rho$, $V$) with one-dimensional center are classified 
in \cite[p.141]{sato-kimura}.  
By differentiating the triplets in the classification \cite[p.141]{sato-kimura}, 
we obtain a classification of isomorphism classes of reduced irreducible (infinitesimal) PVs. 
From this classification 
we choose all reduced 
irreducible PVs of type IFPS,  
then we can obtain a classification of isomorphism classes of reduced irreducible PVs of type IFPS.  
The result is the following:  
\renewcommand{\theenumi}{\alph{enumi}}
\renewcommand{\labelenumi}{(\theenumi)}
\begin{enumerate}
\item $(\gl(1) \oplus \sll(2), \ \Lambda_1 \otimes 3\Lambda_1, \ V(1)\otimes V(4))$
\item $(\gl(1) \oplus \sll(3) \oplus \sll(2), \ \Lambda_1 \otimes 2\Lambda_1 \otimes \Lambda_1, \ V(1) \otimes V(6) \otimes V(2))$
\item $(\gl(1) \oplus \sll(5) \oplus \sll(4), \ \Lambda_1 \otimes \Lambda_2 \otimes \Lambda_1, \ V(1) \otimes V(10)\otimes V(4))$
\end{enumerate}
Here $n \Lambda_1$ is the $n$-th symmetric product of $\Lambda_1$, and $\Lambda_2$ is the second exterior product of $\Lambda_1$.
 
Now let us consider the triplet of the form  
\[(\h, \tau, W) = (\gl(1) \oplus \sll(a) \oplus \sll(a-1), \ \Lambda_1 \otimes \Lambda \otimes \Lambda_1, \ V(1) \otimes V(2a) \otimes V(a-1)),
\] where $a$ is a positive integer and $\Lambda$ is an irreducible representation of $\sll(a)$ with degree $2a$.  
Then the above triplets (a)-(c) are equal to the triplets $(\h, \tau, W)$ for $a$ $=$ 2, 3, or 5 where $\Lambda$ is the 
corresponding irreducible representation.  Note that $\sll(1) = 0$.@

In order to complete our classification, it is sufficient to determine all triplets castling equivalent to the above triplets (a)-(c).  
For this purpose we shall prove Proposition \ref{prop of cast}. 
In the following 
we assume that $m_i$ is a positive natural number.  
\begin{prop}\label{prop of cast} 
Let $a$ be $2$, $3$, or $5$.  
Let $(\g, \rho, V)$  be an arbitrary triplet. 
Then $(\g, \rho, V)$ is castling equivalent to the above triplet $(\h, \tau, W)$  
if and only if $(\g, \rho, V)$ is isomorphic to a triplet of the form 
\begin{eqnarray*}
(\sharp) \quad (\gl(1) \oplus \sll(a) \oplus \sll(m_1) \oplus \cdots \oplus \sll(m_k),
\Lambda_1 \otimes \Lambda \otimes \underbrace{\Lambda_1 \otimes \cdots \otimes \Lambda_1}_k, \\
V(1) \otimes V(2a) \otimes V(m_1) \otimes \cdots \otimes V(m_k)) \ \ (k \geq 1) 
\end{eqnarray*}
which satisfies   
\[(\ast \ast) \ \ a^2 + m_1^2 + \cdots + m_k^2 -k - 2 a m_1 m_2 \cdots m_k = 0.\]
\end{prop}

\pf  
Note that for the triplet $(\h, \tau, W)$,     
we have $\dim$ $\gl(1) \oplus \sll(a) \oplus \sll(a-1)$ $-$ $\dim$ $V(1) \otimes V(2a) \otimes V(a-1)$ $=0$.  
By replacing $a-1$ with $m_1$, we obtain the equality $(\ast \ast)$ $a^2 +m_1^2  -1 -2a m_1=0$. \\[-4mm]

We claim that any triplet castling equivalent to $(\h, \tau ,W)$ is isomorphic to a triplet of the form 
$(\sharp)$ satisfying $(\ast \ast)$. 
Then first note that $(\h, \tau, W)$ itself satisfies the condition $(\sharp)$ and $(\ast \ast)$. 
In the following we prove our claim by induction on the number of sc-transformations.  
 
Now from the hypothesis of induction, 
we consider a triplet $(\g, \rho, V)$ satisfying the condition $(\sharp)$ and $(\ast \ast)$, i.e., 
$(\g, \rho, V)$ is isomorphic to the triplet 
\begin{eqnarray*} 
&&\hspace{-0.8cm} (\h', \tau', W') = (\gl(1) \oplus \sll(a) \oplus \sll(m_1) \oplus \cdots \oplus \sll(m_k), \\ 
&& \hspace{1.5cm} \Lambda_1 \otimes \Lambda \otimes \Lambda_1 \otimes \cdots \otimes \Lambda_1, 
 V(1) \otimes V(2a) \otimes V(m_1) \otimes \cdots V(m_k))
\end{eqnarray*}
and satisfies $(\ast \ast)$.   
Let $(\g', \rho', V')$ be a castling transform of $(\g, \rho, V)$. 
Then it suffices to prove that $(\g', \rho', V')$  is isomorphic to a triplet satisfying $(\sharp)$ and $(\ast \ast)$ again. 
Note that  $(\g', \rho', V')$ is also a castling transform of $(\h', \tau', W')$.  
Then by Proposition \ref{prop iso and cast}, there exists an sc-transformation $(\tilde{\h'}, \tilde{\tau'}, \widetilde{W'})$ of 
$(\h', \tau', W')$ such that $(\tilde{\h'}, \tilde{\tau'}, \widetilde{W'})$ $\cong$ $(\g', \rho', V')$. 
Concerning $(\h', \tau', W')$ there are at most $k$ ways of sc-transformations, however by the symmetry of $m_1, \cdots, m_k$  
we may assume   
that $(\tilde{\h'}, \tilde{\tau'}, \widetilde{W'})$ is given by the triplet
\begin{eqnarray*}
\lefteqn{(\gl(1) \oplus \sll(a) \oplus \sll(m_1) \oplus \cdots \oplus \sll(m_{k-1}) \oplus \sll(2 a m_1 \cdots m_{k-1} - m_k), }\hspace{11cm} \\
\lefteqn{\Lambda_1^* \otimes \Lambda^* \otimes \Lambda_1^* \otimes \cdots \otimes \Lambda_1^* \otimes \Lambda_1,} \hspace{11cm}\\ 
\lefteqn{V(1)^* \otimes V(2a)^* \otimes V(m_1)^* \otimes \cdots \otimes V(m_{k-1})^* \otimes V(2 a m_1 \cdots m_{k-1} - m_k)).} \hspace{11cm} 
\end{eqnarray*}
Then we have
\begin{eqnarray*} 
\lefteqn{ a^2 + m_1^2 + \cdots + m_{k-1}^2 + (2 a m_1 \cdots m_{k-1} - m_k)^2 - k}\hspace{0.6cm} \\
& & {} - 2 a m_1 m_2 \cdots m_{k-1} (2 a m_1 \cdots m_{k-1} - m_k) \\
&=&  a^2 + m_1^2 + \cdots + m_k^2 -k - 2 a m_1 m_2 \cdots m_k \\
&=& 0.
\end{eqnarray*}
Hence $(\tilde{\h'}, \tilde{\tau'}, \widetilde{W'})$ satisfies the equality $(\ast \ast)$.  
Moreover from the claims (1), (4), (5) of Proposition \ref{prop of isom},     
$(\tilde{\h'}, \tilde{\tau'}, \widetilde{W'})$ is isomorphic to a triplet of the form $(\sharp)$.    
Since $(\h, \tau, W)$ satisfies the conditions $(\sharp)$ and $(\ast \ast)$,  
by induction any triplet castling equivalent to $(\h, \tau, W)$ satisfies the conditions $(\sharp)$ and $(\ast \ast)$ again. \\[-1mm] 

Next we prove the converse.  Let $(\g, \rho, V)$ be a triplet isomorphic to the following triplet     
\begin{eqnarray*}
&& \quad (\g', \rho', V')  \\
&& = (\gl(1) \oplus \sll(a) \oplus \sll(m_1) \oplus \cdots \oplus \sll(m_k), 
\Lambda_1 \otimes \Lambda \otimes \Lambda_1 \otimes \cdots \otimes \Lambda_1, \\
&& \quad \ V(1) \otimes V(2a) \otimes V(m_1) \otimes \cdots \otimes V(m_k)) 
\end{eqnarray*}
satisfying the equality $(\ast \ast)$. 
We may assume $2 \leq m_1 \leq m_2 \leq \cdots \leq m_k$ without loss of generality.   
We show that this triplet $(\g', \rho', V')$ is castling equivalent to $(\h, \tau, W)$. 

If $k=1$, the equality $(\ast \ast)$ becomes $a^2 +m_1^2  -1 -2a m_1=0$. The left side of the equality is equal to $(a-m_1)^2 -1$. 
Hence the solution of the equality $(\ast \ast)$ is $m_1 = a \pm 1$. The corresponding triplets are given by 
$(\gl(1) \oplus \sll(a) \oplus \sll(a-1), \ \Lambda_1 \otimes \Lambda \otimes \Lambda_1, \ V(1) \otimes V(2a) \otimes V(a-1))$
and $(\gl(1) \oplus \sll(a) \oplus \sll(a+1), \ \Lambda_1 \otimes \Lambda \otimes \Lambda_1, \ V(1) \otimes V(2a) \otimes V(a+1))$. 
The former is equal to $(\h, \tau, W)$, and the latter is a castling transform of $(\h, \tau, W)$ by claims (1), (4), (5) of Proposition 
\ref{prop of isom}.  Hence $(\g', \rho', V')$ is castling equivalent to $(\h, \tau, W)$. 

Next we investigate the case $k \geq 2$. 
The following lemma plays a crucial role.  The proof of this lemma will be given below. 
\begin{lem}\label{final lem}
Let $k$ be a natural number  such that $k \geq 2$, and let $a$ and $m_i$ $(1 \leq i \leq k)$ be integers satisfying $2 \leq a \leq 5$ and 
$2 \leq m_1 \leq m_2 \leq \cdots \leq m_k$. 
Suppose that the equality 
$(\ast \ast)$ $a^2 + m_1^2 + \cdots + m_k^2 -k - 2 a m_1 m_2 \cdots m_k = 0$ holds.  
Then we have $0 < 2a m_1 m_2 \cdots m_{k-1} -m_k < m_k$. 
\end{lem}

Now by an sc-transformation of the given triplet $(\g', \rho', V')$,  
we can obtain   
the triplet\\[-4mm]
\begin{eqnarray*}
&& \hspace{-3mm}(\g'', \rho'', V'') \\
&&\hspace{-7mm}= (\gl(1) \oplus \sll(a) \oplus \sll(m_1) \oplus \cdots \sll(m_{k-1}) \oplus \sll(2 a m_1 \cdots m_{k-1} - m_k), \\
&& \hspace{-2.8mm}\Lambda_1^* \otimes \Lambda^* \otimes \Lambda_1^* \otimes \cdots \otimes \Lambda_1^* \otimes \Lambda_1, \\ 
&& \hspace{-2.8mm}V(1)^* \otimes V(2a)^* \otimes V(m_1)^* \otimes \cdots \otimes V(m_{k-1})^* \otimes V(2 a m_1 \cdots m_{k-1} - m_k)). 
\end{eqnarray*}
Then by Lemma \ref{final lem}, we have $\dim V'' < \dim V'$.   
Since $\dim V' < \infty$, by a finite number of sc-transformations,  
we can arrive at a triplet of the same form as that of $(\g', \rho', V')$ with some $m_i =1$, i.e. 
\begin{eqnarray*}
(\gl(1) \oplus \sll(a) \oplus \sll(m_1) \oplus \cdots \oplus \sll(1) \oplus \cdots \oplus \sll(m_k),\\
\Lambda_1^{(*)} \otimes \Lambda^{(*)} \otimes \Lambda_1^{(*)} \otimes \cdots \otimes \Lambda_1^{(*)} \otimes \cdots \otimes \Lambda_1^{(*)}, \\
V(1)^{(*)} \otimes V(2a)^{(*)} \otimes V(m_1)^{(*)} \otimes \cdots \otimes V(1)^{(*)} \otimes \cdots \otimes V(m_k)^{(*)}).
\end{eqnarray*}
Thus the number of $\sll(m_i)$-components decreases to $k-1$. Since $k$ is finite, by a finite number of sc-transformations we 
arrive at a triplet isomorphic to  
$(\gl(1) \oplus \sll(a) \oplus \sll(m_1), \ \Lambda_1 \otimes \Lambda \otimes \Lambda_1, \ V(1) \otimes V(2a) \otimes V(m_1))$.    
Since sc-transformations  and a last isomorphism of triplets keep 
the equality $(\ast \ast)$, this triplet also satisfies $(\ast \ast)$. 
Hence by the consideration in the case $k=1$, $(\g', \rho', V')$ is castling equivalent to $(\h, \tau, W)$. 
Therefore ($\g$, $\rho$, $V$) is also castling equivalent to $(\h, \tau, W)$, 
which proves 
Proposition \ref{prop of cast}.  \hfill $\Box$ \\[-1mm]

{\it Proof of Lemma} \ref{final lem}. \ 
Let $k$ be a natural number satisfying $k \geq 2$, and  
let $a$ and  $m_i$ $(1 \leq i \leq k)$ be arbitrary integers satisfying 
the inequalities $2 \leq a \leq 5$,  
$2 \leq m_1 \leq m_2 \leq \cdots \leq m_k$, and 
the equality 
\[
\leqno{(\ast \ast)} \hspace{2cm}  a^2 + m_1^2 + \cdots + m_k^2 -k - 2 a m_1 m_2 \cdots m_k = 0.  
\]

First we show that $0 < 2 a m_1 \cdots m_{k-1} - m_k$. 
We assume that $2 a m_1$ $\cdots m_{k-1}$ $\leq$ $m_k$. Then \ $2 a m_1 \cdots m_{k-1} m_k \leq m_k^2$. 
Hence  
\[a^2 -k + m_1^2 + \cdots + m_k^2 \leq m_k^2.\] 
Therefore 
\[0 \geq a^2 -k + m_1^2 + \cdots + m_{k-1}^2 \geq  a^2 - k + 4(k-1) > 0.\] 
This is a contradiction. Hence we have $0 < 2 a m_1 \cdots m_{k-1} - m_k$.  

Secondly we show that $2 a m_1 m_2 \cdots m_{k-1} -m_k < m_k$. 
Now assume that 
$m_k \leq 3a$. Then since $k \geq 2$ and  $a \geq 2$, we have 
$m_i \leq 2 a 2^{k-1} - 2$. Thus we can use  
Lemma 2 in \cite[p.42]{sato-kimura}, and we have 
\[m_1^2 + \cdots + m_k^2 - 2 a m_1 m_2 \cdots m_k \leq 2^2 k - 2a 2^k.\] 
Substituting this inequality into  $(\ast \ast)$ we have 
\begin{eqnarray*}
0 &=& a^2 + m_1^2 + \cdots + m_k^2 -k - 2 a m_1 m_2 \cdots m_k \\
&\leq& a^2 - k +  2^2 k - 2a 2^k \\
&=& a^2 - 2^{k+1}a + 3 k.  
\end{eqnarray*}
The last expression is negative if $k \geq 2$ and $4- \sqrt{\mathstrut 10} < a < 4 + \sqrt{\mathstrut 10}$, 
thus especially for $2 \leq a \leq 5$. 
This is a contradiction. 
It follows that $m_k > 3a$. 
To prove the inequality $2 a m_1 m_2 \cdots m_{k-1} -m_k < m_k$  we consider two cases: (i) $k=2$ and 
(ii) $k \geq 3$. 

\noindent
(i) First we show that $m_2 > m_1 + a$.  
Now suppose that $m_2 \leq m_1 + a$. 
Then by $(\ast \ast)$ 
\begin{eqnarray*}
0 &=& a^2 + m_1^2 + m_2^2 -2 - 2 a m_1 m_2 \\
&\leq& a^2 + m_1^2 + m_2^2 -2 - 2 a (m_2 - a) m_2 \\
&\leq& a^2 - 2 + (-2a + 2)m_2^2 + 2a^2 m_2 \\
&\leq& a^2 -2 -a m_2^2 + 2a^2 m_2 \\
&=&  a^2 -2 -a m_2(m_2 -2 a). 
\end{eqnarray*}
Since $m_2 - 2a > a$, this expression is less than 
\begin{eqnarray*}
 a^2 - 2 -a^2m_2 < - a^2 - 2 < 0.  
\end{eqnarray*}
This is a contradiction, and 
therefore $m_2 > m_1 +a$. 
Next we show the inequality $2 a m_1 -m_2 < m_2$, which is equivalent to $a m_1 < m_2$. 
Now suppose that $a m_1 \geq m_2$. Then by $(\ast \ast)$ we have 
\begin{eqnarray*}
0 &=& a^2 + m_1^2 + m_2^2 -2 - 2 a m_1 m_2 \\
&\leq& a^2 -2 + m_1^2 + m_2^2 - 2 m_2^2 \\
&=& a^2 -2 + m_1^2 -  m_2^2. 
\end{eqnarray*}
Since $m_2 > m_1 + a$, this expression is less than  
\begin{eqnarray*}
&& a^2 -2 + m_1^2 -  (m_1 + a)^2 \\
&=& -2 -2 a m_1 < 0. 
\end{eqnarray*}
This is a contradiction, and 
therefore we have $2 a m_1 -m_2 < m_2$. 

\noindent
(ii) Consider the case $k \geq 3$. 
First we prove $m_k > \frac{a}{2} (k-1) m_{k-1}$.    
Now suppose that $m_k \leq \frac{a}{2} (k-1) m_{k-1}$. 
Then by $(\ast \ast)$
\begin{eqnarray*}
0 &=& a^2 + m_1^2 + \cdots + m_k^2 -k - 2 a m_1 m_2 \cdots m_k \\
&\leq& a^2 - k + m_1^2 + \cdots + m_k^2 - \dfrac{4}{k-1} m_1 \cdots m_{k-2} m_k^2.
\end{eqnarray*}
Here $\dfrac{4}{k-1} m_1 \cdots m_{k-2} \geq \dfrac{4}{k-1}2^{k-2} = \dfrac{2^k}{k-1} \geq k+1$. 
Hence 
\begin{eqnarray*}
&& a^2 - k + m_1^2 + \cdots + m_k^2 - \dfrac{4}{k-1} m_1 \cdots m_{k-2} m_k^2 \\ 
&\leq& a^2 -k +m_1^2 + \cdots + m_k^2 - (k+1)m_k^2 \\ 
&\leq& a^2 - k - m_k^2.   
\end{eqnarray*}
Since $m_k > 3a > a$, we have $a^2 - k - m_k^2 < a^2 - k -a^2 < 0$.   
This is a contradiction. Hence $m_k > \frac{a}{2} (k-1) m_{k-1}$. 

Next assume that $2 a m_1 \cdots m_{k-1} -m_k$ $\geq$ $m_k$, then $2 a m_1 \cdots m_{k} \geq 2 m_k^2$. 
By applying this inequality to $(\ast \ast)$, we have  
\begin{eqnarray*}
0 &=& a^2 + m_1^2 + \cdots + m_k^2 -k - 2 a m_1 m_2 \cdots m_k \\
&\leq& a^2 - k + m_1^2 + \cdots m_{k-1}^2 - m_k^2. 
\end{eqnarray*}
By $m_k > \frac{a}{2} (k-1) m_{k-1}$,  the last expression is less than  
\begin{eqnarray*}
&& a^2 - k + m_1^2 + \cdots m_{k-1}^2 - \frac{a^2}{4} (k-1)^2 m_{k-1}^2 \\ 
&\leq& a^2 -k - (k-1) (\frac{a^2}{4} (k-1) - 1)m_{k-1}^2 \\
&\leq& a^2 -k - ((k-1)( a^2 (k-1) - 4). 
\end{eqnarray*} 
Finally by $k \geq 3$, we have  
\[0 \leq a^2 -k - 2(2a^2  - 4) = -3a^2 - k + 8 < 0. \]
This is a contradiction, and  
therefore $2 a m_1 m_2 \cdots m_{k-1} -m_k < m_k$. 
\hfill $\Box$
\\[-2mm]

\remark 
Even if $a$ and $m_i$ are real numbers satisfying $2 \leq a \leq 5$ and 
$2 \leq m_1 \leq m_2 \leq \cdots \leq m_k$, we have the same assertion as that of Lemma \ref{final lem} \\

Now we claim that a triplet $(\g, \rho, V)$ is isomorphic to an irreducible PV of type IFPS if and only if 
$(\g, \rho, V)$ is castling equivalent to one of the triplets (a)-(c).  
To see this, first note that a triplet $(\g, \rho, V)$ is isomorphic to an irreducible PV of type IFPS if and only if 
$(\g, \rho, V)$ is castling equivalent to a reduced irreducible PV of type IFPS. 
This follows from the fact that sc-transformations preserve being an irreducible PV of type IFPS and from Corollary \ref{cast of reduced triplet}.   
Furthermore any reduced irreducible PV of type IFPS is isomorphic to one of the triplets (a)-(c). 
Therefore we obtain our claim.  

Thus by applying Proposition \ref{prop of cast} to the case that $(\h, \tau, W)$ is equal to one of the triplets (a)-(c),  
we have the following:    
\begin{cor}\label{classification of irr PV of type IFPS}
Let $k$ and $m_i$ be positive natural numbers. 
A triplet $(\g, \rho, V)$ is isomorphic to an irreducible PV of type IFPS if and only if 
it is isomorphic to a triplet of the form 
\begin{eqnarray*}
&&(\gl(1) \oplus \sll(a) \oplus \sll(m_1) \oplus \cdots \oplus \sll(m_k), \\
&&\Lambda_1 \otimes \Lambda \otimes \underbrace{\Lambda_1 \otimes \cdots \otimes \Lambda_1}_k, \\
&&V(1) \otimes V(2a) \otimes V(m_1) \otimes \cdots \otimes V(m_k))
\end{eqnarray*}
which satisfies the equality 
\[(\ast \ast) \quad  a^2 + m_1^2 + \cdots + m_k^2 -k - 2 a m_1 m_2 \cdots m_k = 0, \] 
where $a$ $=$ $2$, $3$, or $5$, and  
\[\hspace{-2cm}\Lambda = 
\left\{
\begin{array}{ll}
3\Lambda_1 & (a=2) \\
2\Lambda_1 & (a=3) \\
\Lambda_2 & (a=5) \\
\end{array}.
\right.\]
\end{cor}
\vspace{2mm}

\remark 
If a triplet ($\g, \rho, V$) is a PV of type IFPS, then $\dim \g = \dim V$. This equality gives $(\ast \ast)$.   
Indeed 
\[\dim \gl(1) \oplus \sll(a) \oplus \sll(m_1) \oplus \cdots \oplus \sll(m_k) = a^2 + m_1^2 + \cdots + m_k^2 -k,\]
and
\[\hspace{-1.2cm}\dim V(1) \otimes V(2a) \otimes V(m_1) \otimes \cdots \otimes V(m_k) = 2 a m_1 m_2 \cdots m_k.\]
\\[-2mm]

By combining Corollaries \ref{IFPS and PV}, \ref{classification of irr PV of type IFPS} and the claim (3) in Proposition \ref{prop of isom},  
we complete the proof of Theorem \ref{main thm}. \\ 

{\it Remark 1.} \ 
Let $a$ be a natural number such that $a \geq 2$. 
We denote a solution of $(\ast \ast)$ by $(a; m_1, \ldots, m_k)$, which geometrically corresponds to a complex IFPS on 
a Lie group  
$SL(a) \times SL(m_1) \times \cdots \times SL(m_k)$.  
We say that a solution $(a; m_1, \ldots, m_k)$ is essential if $m_i \neq 1$ for all $i$. 
Then for any $k$ we can see that the equation $(\ast \ast)$  
has an essential  solution  by induction.     
First if $k=1$, the equation $(\ast \ast)$ has the solutions $(a; a \pm 1)$. Next we suppose that the equation $(\ast \ast)$ 
has an essential solution $(a; m_1, \ldots, m_k)$.  We may assume that it satisfies $2 \leq m_1 \leq \cdots \leq m_k$.  
When we add $\sll(1)$ to the tail of the Lie algebra $\sll(a) \oplus \sll(m_1) \oplus \cdots \oplus \sll(m_k)$, 
there are $(k+1)$-ways of sc-transformations of this Lie algebra.   
Namely an sc-transformation  of  
$(a; m_1, \ldots, m_k)$ at the $i$-th position for $1 \leq i \leq k$  is given by 
\[(a; m_1, \ldots, m_{i-1}, 2 a m_1 \cdots m_{i-1} m_{i+1} \cdots m_k -m_i, m_{i+1}, \ldots, m_k), \]
and for $i=k+1$ is given by 
\[(a; m_1, \ldots, m_k, 2 a m_1 \cdots m_{k} - 1).\]  
We note that castling transformation of $(a; m_1)$ at the first position is given by $(a; 2 a - m_1)$. 
We have already seen that any sc-transformation of $($$a; m_1, \ldots$, $m_k$$)$    
gives again a solution of $(\ast \ast)$.  
Since we can easily verify that $2 a m_1$ $\cdots m_{k} - 1 > m_k \geq 2$ for the case $i=k+1$, there exists an  
essential solution for any $k \geq 1$. 

For example let us consider the solution $(2; 3)$, which corresponds to a complex IFPS on 
$SL(2) \times SL(3)$.  If we castling transform $(2; 3)$ at the second position,  
we obtain $(2; 3, 11)$. Likewise from ($2$; $3$, $11$) we obtain ($2$; $3, 11, 131$).    

Furthermore even if we fix the number of $\sll$-components $k$ ($\geq 2$), the equation $(\ast \ast)$ has an infinite number 
of solutions as follows:
As we have already showed that there exists at least one non-trivial solution ($a$; $m_1$, $\ldots$, $m_k$). Then by an sc-transformation of 
$(a; m_1, \ldots, m_k)$ at 
the $i$-th position $(1 \leq i < k)$, 
we obtain 
\[(a; m_1, \ldots, m_{i-1}, 2 a m_1 \cdots m_{i-1} m_{i+1} \cdots m_k -m_i, m_{i+1}, \ldots, m_k).\] 
Then this gives a new solution of $(\ast \ast)$ since $2 a m_1 \cdots m_{i-1} m_{i+1} \cdots m_k -m_i > m_k$. 
For example if we castling transform $(2; 3, 11)$ at the first position, we obtain $(2; 41, 11)$. Likewise from $(2; 11, 41)$ 
we obtain $(2; 153, 41)$. 
\\[-1mm]

{\it Remark 2.} \ 
In Theorem \ref{main thm} we have repetition of Lie algebras admitting an IFPS between the cases $a = 2$ and $a = 3$. 
For example $(2; 3)$ and $(3; 2)$, $(2; 3, 11)$ and $(3; 2, 11)$, etc. 
To exclude the repetition, we need the following extra condition:  
for any solution $(3; m_1, \ldots, m_k)$ we require $m_i \neq 2$ for all $i$.    
To see this,  
first we see that by the extra condition we can avoid the repetition. Indeed we can verify that the set of all solutions of $(\ast \ast)$ 
with $a=2$, the one with $a = 3$, the one with $a = 5$ are disjoint. 

To begin with the whole solutions of the cases $a =2$ and $a=3$ obviously does not 
intersect because of the extra condition. 

Next we verify that the whole solutions of the cases $a=2$ and $a =5$ does not intersect as follows:   
By the proof of Proposition \ref{prop of cast}, any solution $(2; m_1, \ldots, m_k)$ of $(\ast \ast)$ is castling equivalent to $(2; 1)$.  
Let us express a solution obtained by a 
castling transformation from $(2; m_1, \ldots, m_k)$ at the $i$-th position as  
($2$; $m_1$, $\cdots$, $m_{i-1}$, $m_i'$, $m_{i+1}$, $\cdots$, $m_k$). 
Then for the case $i = k$ we obtain  
a solution $(2; m_1, \ldots, m_{k-1}, m_k')$ such that $m_k' < m_k$ by Lemma \ref{final lem}. On the other hand for other cases $i \neq k$,   
we have $m_k < m_i'$ by the above consideration.          

Now let us apply this argument to a concrete case. 
We note that 
all solutions obtained by castling transformations from $(2; 1)$ is $(2; 3)$, and those from $(2; 3)$ are $(2; 1)$ and $(2; 3, 11)$. 
Thus concerning any new solution ($2$; $m_1, \ldots, m_k$) castling equivalent to $(2; 3, 11)$, we have $m_i > 11$ for all $i$. 
Hence concerning any solution ($2$; $m_1, \ldots, m_k$) castling equivalent to $(2; 1)$,  we have $m_i \neq 5$ for all $i$. 
Therefore we  do not have repetition between the cases $a =2$ and $a=5$. Similarly we can verify the other case.  

By a similar argument  we can verify that any solution $(3; m_1, \ldots, m_k)$ contains at most one component $m_i$ such as $m_i =2$, 
and in this case changing positions of $m_i$ ($= 2$) and $3$ gives us again the solution 
($2$; $m_1$, $\cdots$, $m_{i-1}$, $3$, $m_{i+1}$, $\cdots$, $m_k$), which is  belonging to the solutions of the case $a=2$. 
It follows  
that the set of all solutions of $(\ast \ast)$ with the condition that $a = 2$, $3$, or $5$ does not change even if 
we add the extra condition. \\[-3mm]

{\it Remark 3.} \ 
Concerning Theorem \ref{main thm} 
we can also obtain an infinite number of real semisimple Lie groups admitting an irreducible IFPS. 
Indeed the triplets obtained in Corollary \ref{classification of irr PV of type IFPS} have real forms 
\begin{eqnarray*}
&& \hspace{-2cm}(\gl(1, \R) \oplus \sll(a, \R) \oplus \sll(m_1,\R) \oplus \cdots \oplus \sll(m_k, \R), \\
&& \hspace{-2cm}\Lambda_1 \otimes \Lambda \otimes \Lambda_1 \otimes \cdots \otimes \Lambda_1, \\
&& \hspace{-2cm}V(1) \otimes V(2a) \otimes V(m_1) \otimes \cdots \otimes V(m_k)). 
\end{eqnarray*}
Hence by Remark 2 after Theorem \ref{(G, X)-str and (G)-hom}, 
real forms $SL(a, \R) \times SL(m_1,\R) \times \cdots \times SL(m_k, \R)$ satisfying  
($\ast \ast$) $a^2 + m_1^2 + \cdots + m_k^2 -k - 2 a m_1 m_2 \cdots m_k = 0$ 
admit an irreducible IFPS.    
Here recall that in \cite{agaoka}, \cite{urakawa}, \cite{elduque}, real simple Lie groups admitting an IFPS were classified, and  
as a result only $\sll(n+1, \R)$ and $\mathfrak{su}^*(2n)$ ($n \geq 1$) admit an IFPS. 
To obtain this result, Aagoka \cite{agaoka} constructed a reducible IFPS on $\sll(n+1, \R)$ and $\mathfrak{su}^*(2n)$ ($n \geq 1$), 
moreover an irreducible IFPS on $\sll(2, \R)$.   
That is why in Theorem \ref{main thm} the complexification of $\sll(n+1, \R)$ ($n \geq 2$) does not appear.  

\section*{Acknowledgment}
The author expresses his sincere thanks to his advisor,  
Professor Yoshio Agaoka at Hiroshima University, for his direction, encouragement and helpful suggestions about his manuscript.      
The author is also indebted to Professor Hideyuki Ishi for his valuable comments about 
a connection between a certain condition of Agaoka \cite[p.145]{agaoka} about left invariant flat projective structures 
and prehomogeneous vector spaces. 
Thanks are due to Professor Tohru Morimoto for several stimulating conversations about 
geometric structures.  
Finally the author expresses his sincere gratitude to the referee who kindly read through the manuscript and gave him many helpful comments.

\end{document}